\documentclass[11pt,reqno]{amsart}
\usepackage[a4paper,margin=30mm]{geometry}
\usepackage[T1]{fontenc}
\usepackage[utf8]{inputenc}
\usepackage{lmodern}
\usepackage{amsmath,amssymb}
\usepackage{booktabs,graphicx}
\newtheorem{theorem}{Theorem}[section]
\newtheorem{lemma}[theorem]{Lemma}
\newtheorem{proposition}[theorem]{Proposition}
\newtheorem{corollary}[theorem]{Corollary}
\newtheorem{definition}[theorem]{Definition}
\theoremstyle{remark}

\numberwithin{equation}{section}
\usepackage{amssymb,amsmath,latexsym}
\usepackage{amsthm}
\usepackage[shortlabels]{enumitem}
\usepackage{caption}
\usepackage{xcolor}
\usepackage{float}
\usepackage{graphicx}
\usepackage{tikz}
\usepackage{url}

\usepackage[colorlinks=true,citecolor=blue,linkcolor=blue,urlcolor=blue,bookmarks,bookmarksopen,bookmarksdepth=2,backref=page,breaklinks]{hyperref}
\usepackage{bookmark}
\bookmarksetup{
	numbered, 
	open,
}

\renewcommand*{\backref}[1]{}
\renewcommand*{\backrefalt}[4]{%
	\ifcase #1 (Not cited.)%
	\or        (Cited on page~#2.)%
	\else      (Cited on pages~#2.)%
	\fi}

\newcommand{\FF}{\mathbb F}
\newcommand{\F}{\mathbb F}
\newcommand{\QQ}{\mathbb Q}
\newcommand{\ZZ}{\mathbb Z}
\newcommand{\CC}{\mathbb C}

\DeclareMathOperator{\Tr}{Tr}
\DeclareMathOperator{\Kl}{Kl}
\allowdisplaybreaks[1]

\title[Walsh spectrum of the Dobbertin function]{The Walsh spectrum of the Dobbertin APN function}
\date{}
\author{Lukas Kölsch}
\address{Department of Mathematics and Statistics,
University of South Florida, St. Petersburg, Florida, USA}
\email{lukas.koelsch.math@gmail.com}

\begin{document}
\begin{abstract}
APN functions yield the best resistance to differential attacks when used as a building block of a block cipher. APN functions induced by monomials over finite fields have received particular interest both for their simplicity and their unique properties. A key property of APN functions is their linearity, which has been determined for all known APN monomials outside of the Dobbertin APN monomials. In this paper, we close this gap and determine the complete Walsh spectrum and in particular the linearity of the Dobbertin APN monomial, answering a problem that has been open for more than 20 years.  Using known connections, our result also yields the asymptotically largest Sidon sets in elementary abelian $2$-groups constructed so far as well as binary linear codes of minimum distance $5$ that are almost optimal.
\end{abstract}
\maketitle
\section{Introduction}
Vectorial Boolean functions play an important role in the design of symmetric cryptosystems, serving as the building blocks of block ciphers. To resist differential attacks, a vectorial Boolean function should have low differential uniformity.
\begin{definition}[Differential uniformity]
	A function $F \colon \F_{2}^n \rightarrow \F_{2}^n$ has differential uniformity $d$, if 
	\begin{equation*}
		d=\max_{a\in \F_{2}^n\setminus\{0\}, b\in \F_{2}^n} |\{x \colon F(x)+F(x+a)=b\}|.
	\end{equation*}
	A function with differential uniformity $2$ is called almost perfect nonlinear (APN).
\end{definition}
As the differential uniformity is always even, APN functions $F \colon \F_{2}^n \rightarrow \F_{2}^n$ provide optimal resistance against differential attacks. One of the most important properties of vectorial Boolean functions for use in substitution permutation networks (SPNs) is their linearity, which can be defined via the Walsh transform of the function.

\begin{definition}[Walsh transform and linearity]
	Let $F \colon \F_{2}^n \rightarrow \F_{2}^n$ be a function. We define its \emph{Walsh transform}
	\begin{equation*}
		W_F(a,b)=\sum_{x \in \F_2^n}(-1)^{\langle a,F(x)\rangle+\langle b,x \rangle},
	\end{equation*}
    where $\langle\cdot ,\cdot \rangle$ is an arbitrary non-degenerate inner product. 
    We call the multiset that collects the values $W_F(a,b)$ for all $a\in \F_2^n\setminus \{0\}$, $b \in \F_2^n$ the \emph{Walsh spectrum} of $F$ and the multiset that collects $|W_F(a,b)|$ for all  $a\in \F_2^n\setminus \{0\}$, $b \in \F_2^n$ the \emph{extended Walsh spectrum} of $F$.
	The \emph{linearity} of $F$ is then defined as $$L(F)=\max_{a \in \F_2^n\setminus\{0\},b\in \F_2^n}|W_F(a,b)|.$$
\end{definition}
    All infinite families of APN functions known so far fall in one of two groups: Either they are quadratic, i.e., of algebraic degree $2$. In this case, the theory of quadratic forms immediately implies that the Walsh transform of these functions can only take on $0$ or powers of $2$. A recent study of the possible Walsh spectra of quadratic APN functions can be found in~\cite{beneteau2026walsh}. The second group of APN families consists of APN functions induced by monomials defined over finite fields, i.e., one identifies $\F_2^n$ with the additive group of the finite field $\F_{2^n}$ and then considers functions of the form $x \mapsto x^d$. A table of all known infinite families of APN exponents (up to CCZ-equivalence, which does not affect the extended Walsh spectrum) can be found in Table~\ref{tab:apn-power-spectra}. As we can see in Table~\ref{tab:apn-power-spectra}, all Walsh spectra of the known APN monomials are known, with the exception of the Dobbertin APN monomials. The main objective of this paper is to fill in this gap. The Dobbertin APN function is a monomial defined on $\F_{q^5}$ where $q=2^m$ via $F_q(x)= x^{D_q}$ where $D_q=q^4+q^3+q^2+q-1$. 

\begin{table}[htbp]
\centering
\small
\setlength{\tabcolsep}{4pt}
\renewcommand{\arraystretch}{1.4}
\resizebox{\textwidth}{!}{%
\begin{tabular}{@{}lllll@{}}
\toprule
Family
& Exponent $d$
& Conditions
& Walsh values
& References\\
\midrule
Gold
& $2^i+1$
& $\gcd(i,n)=1$
& $\begin{aligned}
   &\{0,\pm2^{(n+1)/2}\},
     &&n\text{ odd},\\
   &\{0,\pm2^{n/2},\pm2^{n/2+1}\},
     &&n\text{ even}
   \end{aligned}$
& \cite{Nyberg1994,Gold1968,Yoshiara2017}\\

Kasami
& $2^{2i}-2^i+1$
& $\gcd(i,n)=1$
& $\begin{aligned}
   &\{0,\pm2^{(n+1)/2}\},
     &&n\text{ odd},\\
   &\{0,\pm2^{n/2},\pm2^{n/2+1}\},
     &&n\text{ even}
   \end{aligned}$
& \cite{JanwaWilson1993,Kasami1971,Yoshiara2017}\\

Welch
& $2^t+3$
& $n=2t+1$
& $\{0,\pm2^{(n+1)/2}\}$
& \cite{Dobbertin1999Welch,CanteautCharpinDobbertin2000Welch}\\

Niho
& $2^t+2^{t/2}-1$
& $n=2t+1$, $t$ even
& $\{0,\pm2^{(n+1)/2}\}$
& \cite{Dobbertin1999Niho,HollmannXiang2001}\\

& $2^t+2^{(3t+1)/2}-1$
& $n=2t+1$, $t$ odd
& $\{0,\pm2^{(n+1)/2}\}$
& \cite{Dobbertin1999Niho,HollmannXiang2001}\\

Inverse
& $2^n-2$
& $n$ odd
& $4\mathbb Z\cap
\bigl[1-2^{n/2+1},\,1+2^{n/2+1}\bigr]$
& \cite{Nyberg1994,LachaudWolfmann1990}\\

Dobbertin
& $\begin{gathered}
   2^{4m}+2^{3m}+2^{2m}\\
   {}+2^m-1
   \end{gathered}$
& $n=5m$
& This work
& \cite{Dobbertin2001}\\
\bottomrule
\end{tabular}%
}
\caption{Known APN power families over $\mathbb F_{2^n}$,
$n\geq3$, up to CCZ-equivalence.}
\label{tab:apn-power-spectra}
\end{table}
    Up until this work, the only known result on the Walsh spectrum of the Dobbertin APN function is that all values are divisible by $q^2$ but not all values are divisible by any larger power of $2$, which was shown in 2000 by Canteaut, Charpin, and Dobbertin~\cite{canteaut2000weight}. Since then, no progress on determining the Walsh spectra of the Dobbertin APN function has been made.

\section{Notation and statement of major results}
Fix $q=2^m$ with $m\geq 1$, and put
$K=\FF_q$, $L=\FF_{q^5}$, and $K_r=\FF_{q^r}$.
For a finite field $F$ of characteristic two, write its canonical additive character as
\[
\psi_F(z)=(-1)^{\Tr_{F/\FF_2}(z)}.
\]
This paper is concerned with the Walsh transform of the Dobbertin APN function $F \colon L \rightarrow L$ defined via $F(x)=x^{D_q}$, where $D_q=q^4+q^3+q^2+q-1$. We denote its Walsh transform as
\begin{equation*}
W_q(a,b)=\sum_{x \in \F_{q^5}}(-1)^{\Tr_{\F_{q^5}/\F_2}(ax^{D_q}+bx)}.
\end{equation*}
Note that a simple substitution shows $W_q(a,b)=W_q(a/(b^{D_q}),1)$ if $b\neq 0$, so it actually suffices to determine $W_q(a,1)$ and $W_q(a,0)$. Moreover, we have $\gcd(D_q,q^5-1)=1$ if $q=2^m$ with $m$ odd and $\gcd(D_q,q^5-1)=3$ if $m$ is even (this holds generally for all APN exponents as shown by Dobbertin, see~\cite[Proposition 165]{Carlet2021}), so we get immediately:
\begin{lemma}
    Let $q=2^m$. We have for all $a\neq 0$ that $$W_q(a,0)=\begin{cases}
        (-1)^{m/2}2^{5m/2}, & m \text{ is even, } a \text{ is a non-cube,}\\
        (-1)^{m/2+1}2^{5m/2+1}, & m \text{ is even, } a \text{ is a cube,}\\
        0, & m \text{ is odd.}
    \end{cases}$$
\end{lemma}
\begin{proof}
    By the consideration above, we have $W_q(a,0)=0$ if $m$ is odd since $x \mapsto x^{D_q}$ permutes $L$. If $m$ is even then $W_q(a,0)=W_F(a,0)$ where $F(x)=x^3$. This cubic sum was determined by Carlitz~\cite{carlitz1979explicit}.
\end{proof}
It thus remains to determine $W_q(a,1)$ for all $a \in L^\times$ (note that clearly $W_q(0,1)=0$ and $W_q(0,0)=q^5$). We thus, for simplicity, define $W_q(a):=W_q(a,1)$ for the rest of this paper. 

We will connect the value of $W_q(a,1)$ to the number of points of a specific elliptic curve.

For $0\leq i<5$, set $a_i=a^{q^i}$ and read subscripts modulo five.
Let $e_j(a)$ be the $j$th elementary symmetric polynomial of
$a_0,\ldots,a_4$, and put
\[
e_1(a)=\Tr_{L/K}(a),\qquad e_5(a)=N_{L/K}(a),\qquad
\omega(a)=\Tr_{L/K}(a)^4+\Tr_{L/K}(a)^2e_4(a)+N_{L/K}(a)e_3(a).
\]
These quantities belong to $K$. The single curve used in the proof is
the projective cubic defined by
\begin{equation}\label{eq:curve}
\mathcal E_a:\qquad Y^2+(\Tr_{L/K}(a)X+N_{L/K}(a))Y=X^3+e_2(a)X^2+\omega(a).
\end{equation}
We will prove that $\mathcal E_a$ defines an elliptic curve when $a\notin\{0,1\}$, see Proposition~\ref{prop:nonsingular}. The main goal is to prove the following theorem:

\begin{theorem}\label{thm:main}
For $a\neq 1$, recall the elliptic curve defined via Eq.~\eqref{eq:curve}.
Then for $a\in L^\times$, we have
\[
W_q(a)= \begin{cases}
    q^2(\#\mathcal E_a(\F_q)-q-1), & a\neq 1, \\
    (-1)^{m+1}q^2(q+1), & a=1.
\end{cases}\]
\end{theorem}
The number of rational points of the elliptic curve can be counted using well known techniques to get explicit results. Define the Kloosterman sum over $K$ as
\[
\Kl_K(k)=\sum_{u\in K^\times}\psi_K(u+k/u).
\]
Then we have:
\begin{theorem}\label{thm:main_explicit}
Let $q=2^m$, $K=\F_q$, and $L=\F_{q^5}$.

For $a\in L\setminus\{0,1\}$ with $\Tr_{L/K}(a)\ne0$, put
\[
d(a)=\frac{e_2(a)}{\Tr_{L/K}(a)^2}
     +\frac{N_{L/K}(a)}{\Tr_{L/K}(a)^3},
\qquad
k(a)=
\frac{N_{L/K}\!\left(
\Tr_{L/K}(a)^2+N_{L/K}(a)/a
\right)}
{\Tr_{L/K}(a)^{12}}.
\]
Then $k(a)\ne0$, and for every $a\in L$,
\[
W_q(a)=
\begin{cases}
0,
& a=0,\\[4pt]
(-1)^{m+1}q^2(q+1),
& a=1,\\[6pt]
\displaystyle
q^2\sum_{u\in K}
\psi_K\!\left(
\frac{u^3+e_2(a)u+e_3(a)}{N_{L/K}(a)}
\right),
& a\ne0,\ \Tr_{L/K}(a)=0,\\[12pt]
q^2\psi_K(d(a))\Kl_K(k(a)),
& a\ne1,\ \Tr_{L/K}(a)\ne0.
\end{cases}
\]
In particular, for $a\in K\setminus\{0,1\}$,
\[
W_q(a)
=q^2\psi_K(a)\Kl_K\!\left(\frac{(a+1)^5}{a}\right).
\]
\end{theorem}
Note that the cubic sum in the case where $\Tr_{L/K}(a)=0$ is evaluated explicitly
in Proposition~\ref{prop:cubic-evaluation}, and the Walsh transform is in this case always $0,\pm q^22^{(m+1)/2}$ if $m$ is odd, or $0,\pm q^2 2^{m/2},\pm q^22^{(m+2)/2}$ if $m$ is even. In particular, using the well known characterization of the spectrum of Kloosterman sums, we can determine the extended Walsh spectrum of the Dobbertin APN function precisely:
\begin{theorem}\label{thm:walsh-inclusion}
Put
\[
\mathcal U_m=
\begin{cases}
\{0,2^{(m+1)/2}\},&m\text{ odd},\\
\{0,2^{m/2},2^{m/2+1}\},&m\text{ even}.
\end{cases}
\]
The values occurring in the extended Walsh spectrum of the Dobbertin APN function defined on $L$ are exactly
\[
\begin{aligned}
&\{q^2(q+1)\}\cup q^2\mathcal U_m\cup\{q^2j:j\in\mathbb{N} \text{ odd},\ |j|\leq2\sqrt q\}.
\end{aligned}
\]
In particular, $L(F)=q^2(q+1)$.
\end{theorem}
\section{A reformulation of the sum and further preparation}

The following simple but key lemma transform the character sum attached to the Dobbertin APN function to a cubic sum.
\begin{lemma}\label{lem:permutation}
For every $a\in L$,
\begin{equation*}
W_q(a)=\sum_{y\in L}\psi_L\bigl(ay^{q^2+q+1}+y^{q+1}\bigr).
\end{equation*}
\end{lemma}
\begin{proof}
    A simple calculation (see~\cite[Lemma 9]{budaghyan2022two}) shows that $(q+1)D_q=(q+1)(q^4+q^3+q^2+q-1) \equiv 2q^2(q^2+q+1) \pmod{q^5-1}$, so 
    \begin{align*}
    W_q(a)&=\sum_{y\in L}\psi_L\bigl(ay^{q^4+q^3+q^2+q-1}+y\bigr)\\
    &=\sum_{z\in L}\psi_L\bigl(az^{2q^2(q^2+q+1)}+z^{q+1}\bigr)\\
    &=\sum_{z\in L}\psi_L\bigl(az^{q^2+q+1}+z^{q+1}\bigr),
    \end{align*}
    using the substitution $y=z^{q+1}$ and the fact that $W_q(a)=W_q(a^{2^i})$ for any $i$.
\end{proof}

We also need a simple lemma.
\begin{lemma} \label{lem:auxiliary}
Let $a\in L$, $a \notin \{0,1\}$ then $\Tr_{L/K}(a)^2\neq N_{L/K}(a)/a_i$ for $0\leq i<5$.
\end{lemma}
\begin{proof}
   Assume $\Tr_{L/K}(a)^2=N_{L/K}(a)/a_i$. Since $\Tr_{L/K}(a)^2 \in K$ and $N_{L/K}(a) \in K$ , this implies $a \in K$, and $a_i=a$ for all $i$. Then, in particular, $\Tr_{L/K}(a)=a$ and $N_{L/K}(a)=a^5$, so $N_{L/K}(a)/a_i=a^4$ and $\Tr_{L/K}(a)^2=a^2$, implying $a^2=a^4$, which only holds if $a \in \{0,1\}$, contradicting our assumption. 
\end{proof}

\section{A lemma for the Dobbertin sum over degree $2$ extension fields}
The goal of this section is to relate $W_q(a)$ with $W_{q^2}(a)$, for any $a\in \F_{q^5}^\times$.
\begin{lemma}\label{lem:apn}
For $a\in\FF_{q^5}^\times$,
\begin{equation*}
W_q(a)^2+W_{q^2}(a)=2q^5.
\end{equation*}
\end{lemma}
\begin{proof}
Let $\tau(x)=x^{q^5}$; then $\tau^2=\mathrm{id}$ on $\F_{q^{10}}$ and $\tau$ fixes exactly $\F_{q^5}$. Set $f(x)=x^{D_{q^2}}$ on $\F_{q^{10}}$. Then $f(x)=x^{D_q}$ when $x \in \F_{q^5}^\times$:
Indeed, if $x\in \F_{q^5}^\times$, the five powers $x^{q^{2i}}$ ($0\leq i<5$)
permute $x,x^q,\ldots,x^{q^4}$.

Now consider (unordered) pairs $\{x,x'\}$ of distinct elements of $\F_{q^{10}}$ which are fixed as sets by $\tau$, we call this set $T$. 
If $x,x' \in \F_{q^5}$, then they are clearly fixed by $\tau$ element-wise, so we get
$q^5(q^5-1)/2$ pairs of that shape. Let us call the set of these pairs $T_1\subseteq T$. All other pairs necessarily satisfy $x'=x^{\tau}$, so we get $(q^{10}-q^5)/2$ additional pairs
$\{x,x'\}$ outside $\F_{q^5}$, we call the set of these pairs $T_2\subseteq T$. Thus there are $q^5(q^5-1)/2+(q^{10}-q^5)/2=q^{10}-q^5$ such pairs in total.
The map
\[
h\colon \{x,x'\}\longmapsto (x+x',f(x)+f(x'))
\]
maps $T$ into $\F_{q^5}^\times\times \F_{q^5}$ since $\tau$ fixes both coordinates. We show that $h$ is injective: if 
$h(\{x,x'\})=(s,t)$, then  $f(x+s)+f(x)=t$ and $f(x'+s)+f(x')=t$, and the APN property of the Dobbertin function
does not allow additional solutions. Since $T$ and
$\F_{q^5}^\times\times \F_{q^5}$ have the same cardinality, this map is thus a bijection.
We have
\begin{align*}
W_q(a)^2
&=\sum_{x,x'\in \F_{q^5}}
\psi_{\F_{q^5}}\left(a(f(x)+f(x'))+x+x'\right)\\
&=q^5+
2\sum_{\{x,x'\}\in T_1}
\psi_{\F_{q^5}}\left(a(f(x)+f(x'))+x+x'\right),
\end{align*}
where the $q^5$ term comes from the case $x=x'$.
We also get
\begin{align*}
W_{q^2}(a)&=
\sum_{x\in \F_{q^5}}\psi_{\F_{q^{10}}}(af(x)+x)
+
\sum_{x\in \F_{q^{10}}\setminus \F_{q^5}}\psi_{\F_{q^{10}}}(af(x)+x)\\
&=q^5+\sum_{x\in \F_{q^{10}}\setminus \F_{q^5}}\psi_{\F_{q^{5}}}(af(x)+x+a(f(x))^\tau+x^\tau)\\
&=q^5+2\sum_{\{x,x'\}\in T_2}
\psi_{\F_{q^5}}\left(a(f(x)+f(x'))+x+x'\right),
\end{align*}
where we use in the first step that $af(x)+x \in \F_{q^{5}}$ if $x\in \F_{q^5}$, so $\psi_{\F_{q^{10}}}(af(x)+x)=1$ as well as the transitivity of the trace.
Hence
\[
W_q(a)^2+W_{q^2}(a)=2q^5+2\sum_{\{x,x'\}\in T}
\psi_{\F_{q^5}}\left(a(f(x)+f(x'))+x+x'\right).
\]
Now we use the bijectivity of the map $h$, so we can replace $T$ by $\F_{q^5}^\times \times \F_{q^5}$. We get
\begin{align*}
W_q(a)^2+W_{q^2}(a)&=2q^5+2\sum_{\substack{s\in \F_{q^5}^\times,\\t \in \F_{q^5}}}
\psi_{\F_{q^5}}\left(at+s\right) \\
&=2q^5+2\sum_{s \in \F_{q^5}^\times}\psi_{\F_{q^5}}(s)\sum_{t \in \F_{q^5}}\psi_{\F_{q^5}}(at) \\
&=2q^5,
\end{align*}
using orthogonality of characters (recall $a\neq 0$).
\end{proof}

\section{Relating the Dobbertin sum to another sum}\label{sec:extensions}
With Lemma~\ref{lem:apn} in mind, we try to analyze the sequence $W_{q^r}(a)$ for $r\geq 1$. The idea of this approach is to leverage information from extension fields of $\F_q$ to determine the actual Walsh transform $W_{q}(a)$ (i.e., the case $r=1$ in the sequence). This is a natural tool since our aim is to connect the exponential sum to a point counting problem for an elliptic  curve, whose rational points over extension fields of $\F_q$ can be studied uniformly for all possible extension fields through its zeta function.

The goal of this section is thus to connect the Walsh transform of the Dobbertin function $W_{q^r}(a)$ over arbitrary extension fields $\F_{q^{5r}}$ where $a \in L=\F_{q^5}$ to an exponential sum of a specific polynomial over $\F_{q^r}^5$. 
Fix $a\in L^\times$ and choose $b\in L\setminus K$ arbitrarily. Define linear forms in five variables by
\[
Y_j(X_0,\ldots,X_4)=\sum_{i=0}^4X_i b^{iq^j}, \text{ for } 0\leq j<5.
\]
Since $[L:K]=5$ is prime, $b$ necessarily has degree five over $K$.
Now define
\begin{equation}\label{eq:fixed-polynomial}
P_a(X_0,\dots,X_4)=\sum_{j=0}^4a_jY_jY_{j+1}Y_{j+2}
+\sum_{j=0}^4Y_jY_{j+1}.
\end{equation}
Recall that $a_j=a^{q^j}$.
The map $x \mapsto x^q$ on the coefficients sends $a_j$ to $a_{j+1}$
and $Y_j$ to $Y_{j+1}$, with indices taken modulo $5$. It cyclically permutes both sums, so every
coefficient of $P_a$ is fixed by $x \mapsto x^q$, implying that $P_a\in K[X]$.
We now consider exponential sums using this (fixed) polynomial over arbitrary extension fields of $K$. So, for $r\geq1$, put $K_r=\F_{q^r}$ and set
\begin{equation}\label{eq:extension-sums}
S_r(a)=\sum_{X\in K_r^5}\psi_{K_r}(P_a(X)).
\end{equation}

Let us now consider $S_r(a)$. It turns out that the behavior depends on if $5$ divides $r$ or not (or, equivalently, if $K_r$ contains $L$ or not).

\begin{proposition} \label{prop:split-iden}
 Let $S_r(a)$ be as in Eq.~\eqref{eq:extension-sums}. Then   
 \begin{equation*}
S_r(a)=\begin{cases}
\sum_{z_0,\ldots,z_4\in K_r}
\psi_{K_r}\left(\sum_{i=0}^4a_i z_i z_{i+1}z_{i+2}
+\sum_{i=0}^4z_i z_{i+1}\right),&5\mid r,\\
W_{q^r}(a),&5\nmid r.
\end{cases}
\end{equation*}
\end{proposition}
\begin{proof}
Let us consider first the case that $5\mid r$. Then $b \in K_r$. The coefficient matrix of
the forms $Y_j$ is thus exactly the Vandermonde matrix $(b^{iq^j})_{j,i}$, which has nonzero determinant. The change of variables $(X_0,\dots,X_4)\mapsto (Y_0,\dots,Y_4)$ is thus invertible and we have 
\begin{equation*} 
    S_r(a)=\sum_{X\in K_r^5}\psi_{K_r}(P_a(X))=\sum_{z_0,\ldots,z_4\in K_r}
\psi_{K_r}\left(\sum_{i=0}^4a_i z_i z_{i+1}z_{i+2}
+\sum_{i=0}^4z_i z_{i+1}\right).
\end{equation*}

If $5\nmid r$, then $\{1,b,\ldots,b^4\}$ is a
basis of $\F_{q^{5r}}$ over $K_r$, so the map $(X_0,\dots,X_4)\mapsto Y_0(X_0,\dots,X_4)$ is a bijection between $K_r^5$ and $\F_{q^{5r}}$. Define $y=Y_0(X_0,\dots,X_4) \in \F_{q^{5r}}$. Choose $\ell\in\{1,2,3,4\}$ such that
$\ell r\equiv1\pmod5$ and set $\rho(y)=y^{q^{\ell r}}$.
The first five powers of $\rho$ permute the five Frobenius maps
$z\mapsto z^{q^{jr}}$ ($0\leq j<5$), and $\rho$ maps
$b$ to $b^q$ and $a$ to $a^q$. Since $\rho$ fixes every $X_i\in K_r$,
we necessarily have $Y_j(X_0,\dots,X_4)=\rho^j(y)$ and $a_j=\rho^j(a)$. Hence
\[
P_a(X_0,\dots,X_4)=\Tr_{\F_{q^{5r}}/K_r}\bigl(ay\rho(y)\rho^2(y)+y\rho(y)\bigr).
\]
Put $Q=q^r$ and $M=\F_{Q^5}$.
By \cite[Lemma~9]{budaghyan2022two},
$\gcd(1+Q^\ell,Q^5-1)=1$, and 
\[
D_Q(1+Q^\ell)
\equiv 2Q^{2\ell}(1+Q^\ell+Q^{2\ell})
\pmod{Q^5-1}.
\]
Thus $x=y^{1+Q^\ell}=y\rho(y)$ permutes $M$.
Then
\[
S_r(a)
=\sum_{x\in M}
  \psi_M\bigl(a^{2Q^{2\ell}}x^{D_Q}+x\bigr)
=W_Q(a^{2Q^{2\ell}})
=W_Q(a)
=W_{q^r}(a),
\]
where we use that $W_Q(c^{2^j})=W_Q(c)$.
\end{proof}

 This means that $S_r(a)$ equals the Dobbertin exponential sum over the extension field if $5\nmid r$, and it becomes a much simpler sum in the case $5\mid r$. In particular, $S_1(a)$ is exactly the exponential sum corresponding to the Dobbertin APN function. We will try to determine $S_1(a)$ by considering the sequence of $S_r(a)$, by using the fact that we can calculate $S_r(a)$ in an elementary way using Proposition~\ref{prop:split-iden} if $5 \mid r$. Then we will investigate the sequence of $S_r(a)$, and through this determine $S_1(a)$ and thus  the Walsh transform of the Dobbertin function.

\subsection{The calculation if $5\mid r$}
We now calculate $S_r(a)$ in the case that $5\mid r$.
Recall that $K_r=\FF_{q^r}$.

\begin{proposition}\label{prop:split-eval}
Let  $5\mid r$, and put $Q=q^r$.
For $a\in L^\times$, let $S_r(a)$ be defined by
Eq.~\eqref{eq:extension-sums}.

If $a\notin\{0,1\}$, then
\[
S_r(a)
=
Q^2\left(
\#\mathcal E_a(K_r)-Q-1
-5\psi_{K_r}\bigl((a_1a_4)^{-1}\bigr)
\right),
\]
where $\mathcal E_a$ is the projective curve defined by
\eqref{eq:curve}.

If $a=1$, then
\[
S_r(a)
=
(-1)^{mr}Q^2(4Q-6).
\]
\end{proposition}
The rest of this subsection is dedicated to proving the proposition.

By Proposition~\ref{prop:split-iden}, we need to calculate the sum
\begin{equation} \label{eq:split-sum}
S_r(a)=\sum_{z_0,\ldots,z_4\in K_r}
\psi_{K_r}\left(\sum_{i=0}^4a_i z_i z_{i+1}z_{i+2}
+\sum_{i=0}^4z_i z_{i+1}\right),
\end{equation}
where $a\in L^\times$ and subscripts are read modulo five.
Recall that $a_i=a^{q^i}$ and retain the notation
\[
\Tr_{L/K}(a)=\sum_{i=0}^4a_i,\qquad N_{L/K}(a)=\prod_{i=0}^4a_i.
\]

Since $\psi_{K_r}(u^2)=\psi_{K_r}(u)$, squaring the argument
in Eq.~\eqref{eq:split-sum} leaves its character unchanged.
The substitution $z_i\mapsto z_i^2$ is a bijection of $K_r^5$.
Consequently,
\[
S_r(a)=
\sum_{z_0,\ldots,z_4\in K_r}
\psi_{K_r}\left(
\sum_{i=0}^4 a_i^2 z_i z_{i+1}z_{i+2}
+\sum_{i=0}^4 z_i z_{i+1}
\right).
\]

Expanding the argument gives
\begin{align*}
&a_0^2z_0z_1z_2+a_1^2z_1z_2z_3
+a_2^2z_2z_3z_4+a_3^2z_3z_4z_0
+a_4^2z_4z_0z_1\\
&\hspace{20mm}
+z_0z_1+z_1z_2+z_2z_3+z_3z_4+z_4z_0.
\end{align*}
Set
\[
A=1+a_2^2z_2+a_3^2z_0.
\]
Grouping the terms involving $z_3,z_4$, the argument becomes
\begin{equation}\label{eq:split-phase}
Az_3z_4
+z_2(1+a_1^2z_1)z_3
+z_0(1+a_4^2z_1)z_4
+z_1(a_0^2z_0z_2+z_0+z_2).
\end{equation}

If $A\ne0$, then character orthogonality gives
\begin{align*}
&\sum_{z_3,z_4\in K_r}
\psi_{K_r}\!\left(
Az_3z_4+z_2(1+a_1^2z_1)z_3+z_0(1+a_4^2z_1)z_4
\right)\\
&\quad=
\sum_{z_4\in K_r}\psi_{K_r}(z_0(1+a_4^2z_1)z_4)
\sum_{z_3\in K_r}\psi_{K_r}((Az_4+z_2(1+a_1^2z_1))z_3)\\
&\quad=
Q\psi_{K_r}\!\left(
\frac{z_0z_2(1+a_1^2z_1)(1+a_4^2z_1)}{A}
\right).
\end{align*}
Indeed, the inner sum vanishes unless $Az_4+z_2(1+a_1^2z_1)=0$,
which forces $z_4=z_2(1+a_1^2z_1)/A$.

Including the remaining term in \eqref{eq:split-phase},
the contribution for fixed $z_0,z_2$ with $A\ne0$ is
\begin{equation}\label{eq:split-after-two-variables}
Q\sum_{z_1\in K_r}
\psi_{K_r}\!\left(
z_1(a_0^2z_0z_2+z_0+z_2)
+\frac{z_0z_2(1+a_1^2z_1)(1+a_4^2z_1)}{A}
\right).
\end{equation}

Let $s\in K_r$ be the unique element satisfying
\[
s^2=\frac{z_0z_2}{A}.
\]
Using
\[
(1+a_1^2z_1)(1+a_4^2z_1)
=1+(a_1+a_4)^2z_1+(a_1a_4)^2z_1^2,
\]
the argument in Eq.~\eqref{eq:split-after-two-variables}
becomes
\[
s^2+
\bigl(a_0^2z_0z_2+z_0+z_2+(a_1+a_4)^2s^2\bigr)z_1
+(a_1a_4)^2s^2z_1^2.
\]

Hence the contribution in Eq.~\eqref{eq:split-after-two-variables} equals
\[
Q\psi_{K_r}(s)
\sum_{z_1\in K_r}
\psi_{K_r}\!\left(
\bigl(
a_0^2z_0z_2+z_0+z_2
+(a_1+a_4)^2s^2+a_1a_4s
\bigr)z_1
\right).
\]
Orthogonality shows that this equals $Q^2\psi_{K_r}(s)$
if
\[
a_0^2z_0z_2+z_0+z_2+(a_1+a_4)^2s^2+a_1a_4s=0,
\]
and equals zero otherwise.

Thus the contribution from $A\ne0$, divided by $Q^2$,
is the number of solutions of
\begin{align}
z_0z_2
&=s^2A=s^2(1+a_2^2z_2+a_3^2z_0),
\label{eq:incidence-product}\\
a_0^2z_0z_2+z_0+z_2+(a_1+a_4)^2s^2+a_1a_4s
&=0,
\label{eq:incidence-line}
\end{align}
multiplied by $\psi_{K_r}(s)$.

We first count these equations without imposing $A\ne0$.
Afterward, we subtract their solutions with $A=0$
and restore the contribution from that case in the
original sum.

\subsubsection*{Counting the solutions for fixed $s$}

Fix $s\in K_r$ and introduce 
\[
u=z_0+a_2^2s^2,\qquad v=z_2+a_3^2s^2.
\]
For fixed $s$, these translations are bijections,
so counting pairs $(u,v)$ is equivalent to counting
pairs $(z_0,z_2)$.

We introduce new variables to simplify this system of equations.
\[
p=1+a_0a_3s,\qquad b=1+a_0a_2s,\qquad
j=1+a_2a_3s,\qquad M=a_1a_4+\Tr_{L/K}(a)^2s.
\]

Then
\begin{align*}
    uv&=z_0z_2+a_3^2s^2z_0+a_2^2s^2z_2+a_2^2a_3^2s^4 \\
    &=s^2(1+a_2^2z_2+a_3^2z_0)+a_3^2s^2z_0+a_2^2s^2z_2+a_2^2a_3^2s^4\\
    &=s^2j^2,
\end{align*}
using Eq.~\eqref{eq:incidence-product} in the second step.

Substituting Eq.~\eqref{eq:incidence-product} into Eq.~\eqref{eq:incidence-line} to eliminate $z_0z_2$, and then expressing the result through the newly introduced variables
gives
\[
a_0^2s^2(1+a_2^2z_2+a_3^2z_0)+z_0+z_2+(a_1+a_4)^2s^2+a_1a_4s=
p^2u+b^2v+\Tr_{L/K}(a)^2s^2+a_1a_4s.
\]
Notice further that $p^2u+b^2v+\Tr_{L/K}(a)^2s^2+a_1a_4s=0$ is equivalent to $p^2u+b^2v=sM$.
Thus counting the solutions to Eq.~\eqref{eq:incidence-product} and Eq.~\eqref{eq:incidence-line} is equivalent to counting the solutions of the system of two equations
\begin{equation}\label{eq:sum-product}
uv=s^2j^2,\text{ and } p^2u+b^2v=sM.
\end{equation}
Let $N_s$ denote their number. We distinguish several cases. 

\textbf{Case $sM\ne0$ and $pb\ne0$.}
Assume Eq.~\eqref{eq:sum-product} holds. Then $p^2u$ and $b^2v$ have sum $sM$
and product $(spbj)^2$. Therefore they are the roots of the quadratic equation
\[
w^2+sMw+(spbj)^2=0.
\]
Each root chosen as $p^2u$ determines the other root
$b^2v=sM+p^2u$, and hence determines exactly one
ordered pair $(u,v)$.

The quadratic equation has $2$ roots if $\Tr_{K_r/\F_2}(pbj/M)=0$ and $0$ roots else (this follows from Hilbert's Theorem 90, see e.g~\cite[Section 14.5.]{Carlet2021} for the solvability of quadratic equations in characteristic $2$ more generally), in other words
\begin{equation}\label{eq:smneq0}
    N_s=1+\psi_{K_r}\left(\frac{pbj}{M}\right).
\end{equation}

\textbf{Case $pb=0$, and $sM\ne0$.}
If exactly one of $p,b$ is zero, the second condition of Eq.~\eqref{eq:sum-product}
fixes one of $u,v$ to a unique value (depending on if $p=0$ or $b=0$). The first condition of Eq.~\eqref{eq:sum-product} then fixes the other one, so $N_s=1$ in all these cases.
If both $p,b$ are zero, the  second condition of Eq.~\eqref{eq:sum-product} is $0=sM$,
which is impossible, so $N_s=0$.

\textbf{Case $sM=0$.}
For the moment, assume $a\neq 1$, the case $a=1$ will be handled separately later.
If $s=0$, we have by definition of the auxiliary variables $p=b=j=1$, and
Eq.~\eqref{eq:sum-product} becomes $uv=0$, $u+v=0$.
Thus $u=v=0$ is the unique solution, and $N_0=1$.

If $M=0$, then $\Tr_{L/K}(a)\ne0$ and $s=a_1a_4/\Tr_{L/K}(a)^2$. Plugging this into the definitions of $p$ and $b$ yields
\[
p=\frac{\Tr_{L/K}(a)^2+N_{L/K}(a)/a_2}{\Tr_{L/K}(a)^2},
\qquad
b=\frac{\Tr_{L/K}(a)^2+N_{L/K}(a)/a_3}{\Tr_{L/K}(a)^2}.
\]
Note that both $p,b$ are nonzero. Indeed, if $\Tr_{L/K}(a)^2+N_{L/K}(a)/a_i=0$, then
$a_i=N_{L/K}(a)/\Tr_{L/K}(a)^2\in K$, forcing $a\in K$. In that case
$\Tr_{L/K}(a)=a$, $N_{L/K}(a)=a^5$, and
\[
\Tr_{L/K}(a)^2+\frac {N_{L/K}(a)}{a_i}=a^2(a+1)^2\ne0
\]
because $a\notin\{0,1\}$.

The second condition of Eq.~\eqref{eq:sum-product} now gives $v=(p^2/b^2)u$.
Substituting this into the first  condition of Eq.~\eqref{eq:sum-product} gives
\[
\left(\frac pb u\right)^2=(sj)^2.
\]
There is precisely one solution $u$ of this equation. Thus $N_s=1$ when $M=0$ as well. 

So in total we have $N_s=1+\psi_{K_r}(pbj/M)$ if $pb\neq 0$ and $sM\neq 0$, $N_s=0$ if $p=b=0$ and $N_s=1$ in all other cases.  $p=0$ happens exactly when $s=1/(a_0a_3)$, and $b=0$ exactly when $s=1/(a_0a_2)$. 
\subsubsection*{Summing the weighted counts}
We now calculate $Q^2\cdot N_s \cdot \psi_{K_r}(s)$. For now, we include the case $A=0$, which we will later subtract. By the determination of $N_s$ just now, we get:

\begin{align*}
    \sum_{s \in K_r}Q^2\psi_{K_r}(s)N_s&=Q^2\left(\sum_{\substack{s \in K_r\\sM\neq 0} }\psi_{K_r}(s)(1+\psi_{K_r}(pbj/M))-\psi_{K_r}(1/(a_0a_3))-\psi_{K_r}(1/(a_0a_2))+\sum_{\substack{s\in K_r\\sM=0}}\psi_{K_r}(s)
\right),
\end{align*}
where we include all cases from the previous calculations. This simplifies by character orthogonality to
\begin{equation*}
\begin{aligned}
    \sum_{s \in K_r}Q^2\psi_{K_r}(s)N_s&=Q^2\left(\sum_{\substack{s \in K_r,\\sM\neq 0} }\psi_{K_r}(s)\psi_{K_r}(pbj/M)-\psi_{K_r}(1/(a_0a_3))-\psi_{K_r}(1/(a_0a_2))\right) \nonumber \\
    &=Q^2\left(\sum_{\substack{s \in K_r,\\sM\neq 0} }\psi_{K_r}(s+pbj/M)-2\psi_{K_r}(1/(a_0a_2))\right) \nonumber\\
    &=Q^2\left(\sum_{\substack{s \in K_r,\\M\neq 0} }(\psi_{K_r}(s+pbj/M))-3\psi_{K_r}(1/(a_0a_2))\right) 
\end{aligned}
\end{equation*}
by recognizing that $\psi_{K_r}(1/(a_0a_2))=\psi_{K_r}(1/(a_0a_3))$ since $a_0a_3=(a_0a_2)^{q^3}$, and similarly, for $s=0$ we have $\psi_{K_r}(s+pbj/M)=\psi_{K_r}(1/(a_1a_4))=\psi_{K_r}(1/(a_0a_2))$.

\begin{equation*}
G(s)=s+\frac{pbj}{M}
=s+\frac{(1+a_0a_3s)(1+a_0a_2s)(1+a_2a_3s)}
        {a_1a_4+\Tr_{L/K}(a)^2s}.
\end{equation*}
Hence
\begin{equation}\label{eq:incidence-total}
\sum_{s\in K_r}\psi_{K_r}(s)N_s
=
\sum_{\substack{s\in K_r\\M\ne0}}\psi_{K_r}(G(s))
-3\psi_{K_r}((a_1a_4)^{-1}).
\end{equation}
However, we counted also the solutions for $A=0$, which we have to subtract now.

\subsubsection*{The case $A=0$}

We must remove the  solutions with $A=0$ and
restore their contribution in the original sum.
If $A=0$, Eq.~\eqref{eq:incidence-product}
gives $z_0z_2=0$. Recalling that $A=1+a_2^2z_2+a_3^2z_0$,  the only possibilities are 
$(z_0,z_2)=(0,a_2^{-2})$ or $(z_0,z_2)=(a_3^{-2},0)$.

Eq.~\eqref{eq:incidence-line} then leads to \[
a_j^{-2}+(a_1+a_4)^2s^2+a_1a_4s=0.
\]
for $j=2$ in the first case and $j=3$ in the second case. Multiplying by $a_j^2$ yields a quadratic relation:
\begin{equation}\label{eq:boundary-polynomial}
1+a_j^2a_1a_4s+a_j^2(a_1+a_4)^2s^2=0.
\end{equation}
Under this relation, we have that

\begin{equation*}
\left(\frac{a_1+a_4}{a_1}s+\frac1{a_1a_j}\right)^2
+\left(\frac{a_1+a_4}{a_1}s+\frac1{a_1a_j}\right)
=s+\frac1{a_1a_j},
\end{equation*}

 since, by Eq.~\eqref{eq:boundary-polynomial}, 
$(a_1+a_4)^2s^2/a_1^2+1/(a_1^2a_j^2)=(a_4/a_1)s$. Since the left hand side has absolute trace 0, we deduce that
\[
\psi_{K_r}(s)=\psi_{K_r}((a_1a_j)^{-1}).
\]

If $a_1\ne a_4$, the number of solutions of Eq.~\eqref{eq:boundary-polynomial} is (again using the standard formula for quadratic equations in characteristic 2)
\[
1+\psi_{K_r}\!\left(
\frac{a_1+a_4}{a_ja_1a_4}
\right)
=
1+\psi_{K_r}\!\left(
\frac1{a_1a_j}+\frac1{a_4a_j}
\right).
\] 
Then
\[
\psi_{K_r}(s)\left(1+\psi_{K_r}\!\left(
\frac1{a_1a_j}+\frac1{a_4a_j}
\right)\right)=\psi_{K_r}((a_1a_j)^{-1})
+\psi_{K_r}((a_4a_j)^{-1}).
\]
The total contribution for both $j=2$ and $j=3$ is thus
\begin{equation} \label{eq:total_c1neqc4}
    \psi_{K_r}((a_1a_2)^{-1})
+\psi_{K_r}((a_4a_2)^{-1})+\psi_{K_r}((a_1a_3)^{-1})
+\psi_{K_r}((a_4a_3)^{-1})=2(\psi_{K_r}((a_1a_2)^{-1})
+\psi_{K_r}((a_1a_4)^{-1})).
\end{equation}
 
If $a_1=a_4$, then $a_i=a$  for all $i$ (since $a \in K$) and Eq.~\eqref{eq:boundary-polynomial}  is linear and has
one solution $s$ with  $\psi_{K_r}(s)=\psi_{K_r}(a^{-2})$.

We now compare these excluded weights with the original
 contribution at $A=0$. In Eq.~\eqref{eq:split-phase}, the
sum over $z_3,z_4$ at $A=0$ equals $Q^2$ precisely when
$z_2(1+a_1^2z_1)=0$ and 
$z_0(1+a_4^2z_1)=0$.

If $a_1\ne a_4$, then $z_0,z_2$ cannot both be nonzero.
There are thus only two possible solutions
$(z_0,z_1,z_2)=(0,a_1^{-2},a_2^{-2})$ and
$(z_0,z_1,z_2)=(a_3^{-2},a_4^{-2},0)$,
recalling that $z_0=z_2=0$ would lead to $A=1$. 
In both cases the character $\psi_{K_r}$ of the term in Eq.~\eqref{eq:split-phase} is $\psi_{K_r}((a_1a_2)^{-1})$ (after applying an appropriate Frobenius automorphism inside the character). 

Thus the net correction is, taking Eq.~\eqref{eq:total_c1neqc4} into account,
$-2\psi_{K_r}((a_1a_4)^{-1})$ for the case $a_1\neq a_4$.

If $a_1=a_4$, then $z_2(1+a_1^2z_1)=0$ and 
$z_0(1+a_4^2z_1)=0$ implies (keeping in mind that now $a_i=a$ for all $i$ and $z_0=z_2=0$ is impossible) $z_1=a^{-2}$ and the condition $A=1+a^2(z_0+z_2)=0$ forces
$z_2=z_0+a^{-2}$. 
Plugging these values into Eq.~\eqref{eq:split-phase} yields a remaining term of 
$z_0^2+a^{-2}z_0+a^{-4}$. Therefore its contribution is
\begin{equation*}
\begin{aligned}
Q^2\sum_{z_0\in K_r}
\psi_{K_r}(z_0^2+a^{-2}z_0+a^{-4})
&=Q^2\psi_{K_r}(a^{-1})
\sum_{z_0\in K_r}\psi_{K_r}((1+a^{-2})z_0)
\nonumber\\
&=
\begin{cases}
0,&a\ne1,\\
Q^3\psi_{K_r}(1),&a=1.
\end{cases}
\end{aligned}
\end{equation*}

Consequently, for every $a\in L^\times$, the entire
boundary correction is accounted for by
\begin{equation}\label{eq:split-boundary-correction}
S_r(a)
=
Q^2\left(\sum_{s\in K_r}\psi_{K_r}(s)N_s
-2\psi_{K_r}((a_1a_4)^{-1})
+\epsilon Q\psi_{K_r}(1)\right),
\end{equation}
where $\epsilon=1$ if $a=1$ and $0$ else.
For $a\notin\{0,1\}$, combining this with
Eq.~\eqref{eq:incidence-total} gives
\begin{equation}\label{eq:split-G}
S_r(a)
=
Q^2\left(\sum_{\substack{s\in K_r,\\a_1a_4+\Tr_{L/K}(a)^2s\ne0}}
\psi_{K_r}\left(s+\frac{(1+a_0a_3s)(1+a_0a_2s)(1+a_2a_3s)}
        {a_1a_4+\Tr_{L/K}(a)^2s}\right)
-5\psi_{K_r}((a_1a_4)^{-1})\right).
\end{equation}
\subsubsection*{Comparing the remaining sum with the curve}

Continue for now to assume $a\notin\{0,1\}$.
We interpret the character sum in Eq.~\eqref{eq:split-G} as the number of points on a curve.

For each $s\in K_r$ with $a_1a_4+\Tr_{L/K}(a)^2s\ne0$, the equation
\begin{equation}\label{eq:split-artin-schreier}
w^2+w
=
s+\frac{(1+a_0a_3s)(1+a_0a_2s)(1+a_2a_3s)}
        {a_1a_4+\Tr_{L/K}(a)^2s}
\end{equation}
has exactly
\[
1+\psi_{K_r}\!\left(
s+\frac{(1+a_0a_3s)(1+a_0a_2s)(1+a_2a_3s)}
        {a_1a_4+\Tr_{L/K}(a)^2s}
\right)
\]
solutions.

We first rewrite this equation as a cubic.
Multiplying Eq.~\eqref{eq:split-artin-schreier} by
$(a_1a_4+\Tr_{L/K}(a)^2s)^2$ gives
\begin{align*}
&\bigl((a_1a_4+\Tr_{L/K}(a)^2s)w\bigr)^2
+(a_1a_4+\Tr_{L/K}(a)^2s)\bigl((a_1a_4+\Tr_{L/K}(a)^2s)w\bigr)\\
&=
s(a_1a_4+\Tr_{L/K}(a)^2s)^2+
(a_1a_4+\Tr_{L/K}(a)^2s)
(1+a_0a_3s)(1+a_0a_2s)(1+a_2a_3s).
\end{align*}

The right side has degree at most four in $s$, with leading term
\[
\Tr_{L/K}(a)^2a_0^2a_2^2a_3^2s^4
=(\Tr_{L/K}(a)a_0a_2a_3s^2)^2.
\]
We remove this square term by setting
\[
v=(a_1a_4+\Tr_{L/K}(a)^2s)w+\Tr_{L/K}(a)a_0a_2a_3s^2.
\]
For each fixed $s$ with $a_1a_4+\Tr_{L/K}(a)^2s\ne0$, this is a
bijection between $w$ and $v$.

Substitution gives
\[
\begin{aligned}
v^2+(a_1a_4+\Tr_{L/K}(a)^2s)v
&=
s(a_1a_4+\Tr_{L/K}(a)^2s)^2\\
&\quad+
(a_1a_4+\Tr_{L/K}(a)^2s)
(1+a_0a_3s)(1+a_0a_2s)(1+a_2a_3s)\\
&\quad+
(\Tr_{L/K}(a)a_0a_2a_3s^2)^2
+(a_1a_4+\Tr_{L/K}(a)^2s)\Tr_{L/K}(a)a_0a_2a_3s^2.
\end{aligned}
\]
Expanding and cancelling the degree-four terms yields
the cubic
\begin{equation}\label{eq:split-cubic-model}
\begin{aligned}
v^2+(a_1a_4+\Tr_{L/K}(a)^2s)v
&=
\left(\Tr_{L/K}(a)^2+\frac {N_{L/K}(a)}{a_1}\right)
\left(\Tr_{L/K}(a)^2+\frac {N_{L/K}(a)}{a_4}\right)s^3\\
&\quad+
\bigl(
\Tr_{L/K}(a)^2(a_0a_2+a_0a_3+a_2a_3)
+N_{L/K}(a)(a_1+a_4)
\bigr)s^2\\
&\quad+
\bigl(
\Tr_{L/K}(a)^2+(a_1a_4)^2
+a_1a_4(a_0a_2+a_0a_3+a_2a_3)
\bigr)s\\
&\quad+a_1a_4,
\end{aligned}
\end{equation}
where we use $N_{L/K}(a)=a_0a_1a_2a_3a_4$.

For $a\neq 1$, the factors $\Tr_{L/K}(a)^2+N_{L/K}(a)/a_i$ are
nonzero by Lemma~\ref{lem:auxiliary}.
In particular, the coefficient of $s^3$
in Eq.~\eqref{eq:split-cubic-model} is nonzero.

Now make the affine change of coordinates
\begin{equation*}
\begin{aligned}
X&=
\left(\Tr_{L/K}(a)^2+\frac {N_{L/K}(a)}{a_1}\right)
\left(\Tr_{L/K}(a)^2+\frac {N_{L/K}(a)}{a_4}\right)s
+a_1a_4\Tr_{L/K}(a)^2+N_{L/K}(a)(a_1+a_4),\\
Y&=
\left(\Tr_{L/K}(a)^2+\frac {N_{L/K}(a)}{a_1}\right)
\left(\Tr_{L/K}(a)^2+\frac {N_{L/K}(a)}{a_4}\right)v\\
&\quad+
\bigl(e_2(a)+\Tr_{L/K}(a)(a_1+a_4)\bigr)X
+\omega(a)+\Tr_{L/K}(a)a_0a_2a_3(a_1a_4)^2.
\end{aligned}
\end{equation*}
Substituting this into Eq.~\eqref{eq:split-cubic-model}, using
the definitions of $e_2(a),e_3(a),e_4(a)$ and
$\omega(a)=\Tr_{L/K}(a)^4+\Tr_{L/K}(a)^2e_4(a)+N_{L/K}(a)e_3(a)$, gives
\begin{equation*}
Y^2+(\Tr_{L/K}(a)^2X+N_{L/K}(a)^2)Y
=X^3+e_2(a)^2X^2+\omega(a)^2.
\end{equation*}

This is the equation obtained by squaring every
coefficient of
\[
\mathcal E_a:\qquad
y^2+(\Tr_{L/K}(a)x+N_{L/K}(a))y=x^3+e_2(a)x^2+\omega(a).
\]
Squaring all coefficients in characteristic $2$ of course preserves the number of points of a curve, hence the projective cubic in 
\eqref{eq:split-cubic-model} has exactly
$\#\mathcal E_a(K_r)$ points (note that both have exactly one point at infinity).

If $a_1a_4+\Tr_{L/K}(a)^2s=0$, 
Eq.~\eqref{eq:split-cubic-model} has a unique solution in $v$.
The number of affine points on $\mathcal E_a(K_r)$ is thus exactly
\begin{align*}
&\sum_{\substack{s\in K_r,\\a_1a_4+\Tr_{L/K}(a)^2s\ne0}}\left(1+\psi_{K_r}\left(
s+\frac{(1+a_0a_3s)(1+a_0a_2s)(1+a_2a_3s)}
        {a_1a_4+\Tr_{L/K}(a)^2s}
\right)\right)+\sum_{\substack{s\in K_r,\\a_1a_4+\Tr_{L/K}(a)^2s=0}}1\\&=Q+\sum_{\substack{s\in K_r,\\a_1a_4+\Tr_{L/K}(a)^2s\ne0}}\psi_{K_r}\left(
s+\frac{(1+a_0a_3s)(1+a_0a_2s)(1+a_2a_3s)}
        {a_1a_4+\Tr_{L/K}(a)^2s}
\right)
\end{align*}

Taking the point at infinity into account, we get
\begin{equation*}
\begin{aligned}
\#\mathcal E_a(K_r)-Q-1
&=
\sum_{\substack{s\in K_r,\\a_1a_4+\Tr_{L/K}(a)^2s\ne0}}
\psi_{K_r}\!\left(
s+\frac{(1+a_0a_3s)(1+a_0a_2s)(1+a_2a_3s)}
        {a_1a_4+\Tr_{L/K}(a)^2s}
\right).
\end{aligned}
\end{equation*}

Combining this identity with Eq.~\eqref{eq:split-G}, we obtain
\begin{equation*}
S_r(a)
=
Q^2\left(
\#\mathcal E_a(K_r)-Q-1
-5\psi_{K_r}((a_1a_4)^{-1})
\right).
\end{equation*}
\subsubsection*{The case $a=1$}

Now suppose $a=1$, so all $a_i=1$ and $\Tr_{L/K}(a)=N_{L/K}(a)=1$.
The conditions in Eq.~\eqref{eq:sum-product} become
\[
\begin{cases}
(z_0+s^2)(z_2+s^2)=s^2(1+s)^2,\\
(1+s)^2(z_0+z_2)=s(1+s).
\end{cases}
\]
If $s=0$, there is one solution the equations reduce to $z_0z_2=0$, $z_0+z_2=0$, so there is only one solution $z_0=z_2=0$.
If $s=1$, the equations reduce to
$(z_0+1)(z_2+1)=0$, so $z_0=1$ or $z_2=1$ with the other one chosen arbitrarily. We have thus in this case $2Q-1$ solutions.
If $s\notin\{0,1\}$, we have (using the previous notation) $sM\neq 0$ and $pb\neq 0$, so by Eq.~\eqref{eq:smneq0} the number of solutions is $1+\psi_{K_r}((1+s)^2)$.

Thus we have in total

\begin{align*}
\sum_{s\in K_r}\psi_{K_r}(s)N_s
&=
1+(2Q-1)\psi_{K_r}(1)
+\sum_{s\in K_r\setminus\{0,1\}}
\left(\psi_{K_r}(s)+\psi_{K_r}(s)\psi_{K_r}((1+s)^2)\right)\\
&=1+(2Q-1)\psi_{K_r}(1)
+\sum_{s\in K_r\setminus\{0,1\}}
\psi_{K_r}(s)+\sum_{s\in K_r\setminus\{0,1\}}\psi_{K_r}(1)\\
&=1+(2Q-1)\psi_{K_r}(1)-(1+\psi_{K_r}(1))+(Q-2)\psi_{K_r}(1) \\
&=(3Q-4)\psi_{K_r}(1),
\end{align*}
using orthogonality.

The case for $A=0$ is dealt with in 
Eq.~\eqref{eq:split-boundary-correction}, so we obtain
\[
S_r(a)
=Q^2\left((3Q-4)\psi_{K_r}(1)-2\psi_{K_r}(1)+Q\psi_{K_r}(1)\right)
=Q^2(4Q-6)\psi_{K_r}(1)=(-1)^{mr}Q^2(4Q-6).
\]

\section{The corresponding elliptic curve}\label{sec:curve}
Let $a\in L^\times$ be fixed. Write (compare Eq.~\eqref{eq:curve})

\[f_a(X,Y)=Y^2+(\Tr_{L/K}(a)X+N_{L/K}(a))Y+X^3+e_2(a)X^2+\omega(a).\]

The corresponding homogeneous polynomial is
\[
F_a(X,Y,Z)=Y^2Z+(\Tr_{L/K}(a)X+N_{L/K}(a)Z)YZ+X^3+e_2(a)X^2Z+\omega(a) Z^3.
\]
\begin{proposition} \label{prop:nonsingular}
    The projective cubic variety $V$ defined by the equation $F_a(X,Y,Z)=0$ is nonsingular for any $a \in L\setminus \{0,1\}$.
\end{proposition}
\begin{proof}
At $Z=0$ it has just the point $O=[0:1:0]$. We have the partial derivatives 
\[\frac{\partial F_a}{\partial X}
=\Tr_{L/K}(a)YZ+X^2, \qquad
\frac{\partial F_a}{\partial Y}
=\Tr_{L/K}(a)XZ+N_{L/K}(a)Z^2, \qquad \frac{\partial F_a}{\partial Z}=Y^2+\Tr_{L/K}(a)XY+e_2(a)X^2+\omega(a) Z^2.
\]
At $O=[0:1:0]$, the partial derivative $\frac{\partial F_a}{\partial Z}$ clearly does not vanish, so $O$ is nonsingular. 
So let us check the affine points (i.e., via scaling, we may choose $Z=1$). 

If $\Tr_{L/K}(a)=0$ then $\frac{\partial F_a}{\partial Y}$ does not vanish because $N_{L/K}(a)\ne0$.

If $\Tr_{L/K}(a)\ne0$ and $\frac{\partial F_a}{\partial X}(X,Y,1)=\frac{\partial F_a}{\partial Y}(X,Y,1)=0$ then $X=N_{L/K}(a)/\Tr_{L/K}(a)$ and $Y=N_{L/K}(a)^2/\Tr_{L/K}(a)^3$. For these choices of $X,Y$ we have 
\[F_a(X,Y,1)=f_a(X,Y)=\frac{N_{L/K}(a)^4}{\Tr_{L/K}(a)^6}+\frac{N_{L/K}(a)^3}{\Tr_{L/K}(a)^3}+e_2(a)\frac{N_{L/K}(a)^2}{\Tr_{L/K}(a)^2}+\omega(a)=\frac{N_{L/K}(a)^4+N_{L/K}(a)^3\Tr_{L/K}(a)^3+e_2(a)N_{L/K}(a)^2\Tr_{L/K}(a)^4+\omega(a) \Tr_{L/K}(a)^6}{\Tr_{L/K}(a)^6}.\]
It is thus enough to show that the numerator does not vanish.
We have
\begin{align*}
N_{L/K}(a)^4+N_{L/K}(a)^3\Tr_{L/K}(a)^3+e_2(a)N_{L/K}(a)^2\Tr_{L/K}(a)^4+\omega(a) \Tr_{L/K}(a)^6&=N_{L/K}(a)^4+N_{L/K}(a)^3\Tr_{L/K}(a)^3+e_2(a)N_{L/K}(a)^2\Tr_{L/K}(a)^4+e_3(a)N_{L/K}(a)\Tr_{L/K}(a)^6+e_4(a)\Tr_{L/K}(a)^8+\Tr_{L/K}(a)^{10}\\
&=\frac1{N_{L/K}(a)}\prod_{i=0}^4(N_{L/K}(a)+\Tr_{L/K}(a)^2a_i),
\end{align*}
where we use $\omega(a)=\Tr_{L/K}(a)^4+\Tr_{L/K}(a)^2e_4(a)+N_{L/K}(a)e_3(a)$, $e_1(a)=\Tr_{L/K}(a)$, $e_5(a)=N_{L/K}(a)$, and $\prod_i a_i=N_{L/K}(a)$. By Lemma~\ref{lem:auxiliary}, $N_{L/K}(a)+\Tr_{L/K}(a)^2a_i\neq 0$ if $a\notin \{0,1\}$, so the term cannot vanish, proving that $V$ has no singular point in these cases.
\end{proof}

A nonsingular projective cubic variety has genus one, and the point
$O$ makes it an elliptic curve over $K$, see~\cite[Proposition~III.3.1(c)]{Silverman2009}, which is exactly the curve given in Eq.~\eqref{eq:curve}:
\[\mathcal E_a:\qquad Y^2+(\Tr_{L/K}(a)X+N_{L/K}(a))Y=X^3+e_2(a)X^2+\omega(a),\]
as long as $a\neq 1$.
So let for now $a \in L\setminus \{0,1\}$. For any finite extension $K_r/K$, write
\begin{equation} \label{eq:curve_count}
    C_r(a)=\#\mathcal E_a(K_r)-|K_r|-1,
\end{equation}

where $\#\mathcal E_a(K_r)$ denotes the number of $K_r$-rational points of $\mathcal E_a$.

For $a=1$ we are not able to use the connection to the elliptic curve as in this case the projective variety is singular. Instead, we simply define directly
\begin{equation}\label{eq:C-exceptional}
C_r(1)=(-1)^{rm+1}(q^r+1).
\end{equation}

\begin{theorem} \label{thm:two-roots}
    For every $a \in L^\times$, there exist
$\alpha,\beta\in\CC$ such that, for every $r\geq1$,
\begin{equation*}
C_r(a)=-\alpha^r-\beta^r,\qquad
\alpha+\beta=-C_1(a)\in\ZZ,\qquad \alpha\beta=q.
\end{equation*}
Consequently,
\begin{equation*}
C_r(a)^2+C_{2r}(a)
=(\alpha^r+\beta^r)^2-\alpha^{2r}-\beta^{2r}
=2(\alpha\beta)^r
=2q^r.
\end{equation*}
\end{theorem}
\begin{proof}
    For $a\neq 1$, this follows from the standard theory of the zeta functions of elliptic curves, see e.g~\cite[Theorem~V.2.3.1(a)]{Silverman2009}. For $a=1$, we can find $\alpha=(-1)^m q$ and $\beta=(-1)^m$ directly.
\end{proof}

\section{Investigating the exponential sum over extension fields}
In this section, we establish a relation for the sums $S_r(a)$ through the application of Dwork's theorem. Outside of Dwork's theorem itself (we refer to both Dwork's original work~\cite{Dwork} as well as the standard reference~\cite{Koblitz1984}), we will explain all necessary (mostly standard) steps for the convenience of potential readers that are not familiar with these techniques. 
Dwork's rationality theorem  (we will use the version found in \cite[Theorem on p. 113]{Koblitz1984}) states that, for an affine hypersurface
$V$ over $K$, its zeta function can be expressed as a rational function, namely, 
\[
Z_V(T)=\exp\left(\sum_{r\geq1}\#V(K_r)\frac{T^r}{r}\right)=\frac{P(T)}{Q(T)},
\]
where $P,Q \in \ZZ[T]$ with $P(0)=Q(0)=1$. Here $\#V(K_r)$ is the number of $K_r$-rational points of $V$. The following proposition is then fairly standard to prove, see e.g.~\cite[Exercise 6, p. 115]{Koblitz1984} for a similar variant.

\begin{proposition} \label{prop:dwork}
There exist finitely many $\gamma_1,\dots,\gamma_k \in \CC \setminus \{0\}$, $m_1,\dots,m_k \in \ZZ$ such that 
    \[
S_r(a)=\sum_{i=1}^k m_i\gamma_i^r,
\]
for every $r\geq 1$.
\end{proposition}
\begin{proof}
We define the hypersurface 
\[
V_a:\qquad w^2+w=P_a(X_0,\ldots,X_4)\]
over $K$, where $P_a$ is defined as in Eq.~\eqref{eq:fixed-polynomial}.
For each $(X_0,\dots,X_4)\in K_r^5$, the equation $w^2+w=P_a(X_0,\ldots,X_4)$ has exactly $1+\psi_{K_r}(P_a(X_0,\dots,X_4))$ solutions, so the total number of $K_r$-rational points on $V_a$ is

\begin{align*}
\#V_a(K_r)
&=\sum_{X_0,\dots,X_4\in K_r}
  \bigl(1+\psi_{K_r}(P_a(X_0,\dots,X_4))\bigr)\\
&=q^{5r}+S_r(a),
\end{align*}
by the definition of $S_r(a)$ in Eq.~\eqref{eq:extension-sums}.

We now apply Dwork's theorem to the zeta function of $V_a$. Write $Z_{V_a}(T)=P(T)/Q(T)$. We factor both $P$ and $Q$ in linear factors of the form $(1-\lambda T)$ with $\lambda \in \CC\setminus\{0\}$. Collecting all these factors of both $P$ and $Q$ gives a product
\[
Z_{V_a}(T)=\prod_\lambda(1-\lambda T)^{-m_\lambda},
\]
with finitely many nonzero $\lambda$. Here $m_\lambda \in \ZZ$ is the multiplicity of the factor
$1-\lambda T$ in $Q$ minus its multiplicity
in $P$. Using the formal power series yields 
$-\log(1-\lambda T)
=\sum_{r\geq1}\lambda^r\frac{T^r}{r},$ so

\begin{align*}
\log Z_{V_a}(T)
&=-\sum_\lambda m_\lambda\log(1-\lambda T)\\
&=\sum_{r\geq1}
  \left(\sum_\lambda m_\lambda\lambda^r\right)
  \frac{T^r}{r}.
\end{align*}
On the other hand, by the definition of the zeta function, 
\begin{align*}
\log Z_{V_a}(T)
&=\sum_{r\geq1}\#V_a(K_r)\frac{T^r}{r}\\
&=\sum_{r\geq1}(q^{5r}+S_r(a))\frac{T^r}{r}.
\end{align*}
Comparing coefficients therefore yields
\[
q^{5r}+S_r(a)=\sum_\lambda m_\lambda\lambda^r.
\]
Of course $q^{5r}$ can itself be written as $m_\lambda\lambda^r$ by picking $\lambda=q^5$, proving the statement.
\end{proof}

We now remove the additional terms in Proposition~\ref{prop:split-eval}
from the extension sums $S_r(a)$. Define
\begin{equation}\label{eq:corrected-sequence}
\mathcal A_r(a)=
\begin{cases}
-S_r(a),
&5\nmid r,\\
-S_r(a)-5q^{2r}\psi_{K_r}\bigl((a_1a_4)^{-1}\bigr),
&5\mid r,\ a\ne1,\\
-S_r(a)+5q^{2r}(q^r-1)(-1)^{mr},
&5\mid r,\ a=1.
\end{cases}
\end{equation}
By Proposition~\ref{prop:split-eval}, Eq.~\eqref{eq:curve_count} and Eq.~\eqref{eq:C-exceptional},
these corrections give 
$\mathcal A_r(a)=-q^{2r}C_r(a)$ if $5\mid r$.
Consequently, Theorem~\ref{thm:two-roots} gives for all $r\geq 1$
\begin{equation}\label{eq:fifth-subsequence}
\mathcal A_{5r}(a)=(q^2\alpha)^{5r}+(q^2\beta)^{5r}
\end{equation}
We will use this subsequence of the sequence of $\mathcal A_r(a)$ to determine $\mathcal A_1(a)$ and thus $S_1(a)$, which is precisely the Dobbertin APN exponential sum.
More precisely, our goal is to prove $\mathcal A_1(a)=-q^2C_1(a)$, which gives
$S_1(a)=q^2C_1(a)$.
The first step is to show that the correction preserves
the finite representation supplied by Dwork's theorem in Proposition~\ref{prop:dwork}.

\begin{proposition}\label{prop:corrected-power-sum}
The sequence $(\mathcal A_r(a))_{r\geq1}$ is integer-valued and there exist finitely many $\gamma_1,\dots,\gamma_k \in \CC \setminus \{0\}$, $m_1,\dots,m_k \in \ZZ$ such that
\begin{equation}\label{eq:a-power-sum}
\mathcal A_r(a)=\sum_{i=1}^km_i\gamma_i^r,
\end{equation}
where additionally
\begin{equation}\label{eq:total-multiplicity}
\sum_{i=1}^km_i=2.
\end{equation}
\end{proposition}

\begin{proof}
The integrality of the sequence follows directly by the definition in Eq.~\eqref{eq:corrected-sequence}.

Proposition~\ref{prop:dwork} gives a finite power-sum representation
for $S_r(a)$, and hence also for $-S_r(a)$, say
\begin{equation} \label{eq:orig_sum}
    -S_r(a)=\sum_{i=1}^lt_i{\gamma}_i^r,
\end{equation}

where $t_i \in \ZZ$.

Put
\[
\chi=\psi_L((a_1a_4)^{-1}).
\]
If $r=5h$, trace transitivity gives
\[
\psi_{K_r}((a_1a_4)^{-1})
=\chi^h=\chi^r,
\]
because $\chi\in\{1,-1\}$ and $h,r$ have the same parity.
Thus the definition of $\mathcal A_r(a)$ in Eq.~\eqref{eq:corrected-sequence} gives
\begin{equation}\label{eq:boundary-correction}
\mathcal A_r(a)+S_r(a)=
\begin{cases}
-5(q^2\chi)^r,
&5\mid r, a\ne1,\\
5
\bigl((q^3(-1)^m)^r-(q^2(-1)^m)^r\bigr),
&5\mid r, a=1, \\
0, & 5 \nmid r.
\end{cases}
\end{equation}

Let $\zeta$ be a primitive complex fifth root of unity. Then, for every
$z\in\CC^\times$,
\begin{equation*}
\sum_{j=0}^4(z\zeta^j)^r
=\begin{cases}
5z^r,
&5\mid r\\
0,
& \text{else.}
\end{cases}
\end{equation*}
Applying this identity with $z=q^2\chi$ in the first case of Eq.~\eqref{eq:boundary-correction}, and using both
$z=q^3(-1)^m$ and $z=q^2(-1)^m$ in the second, shows that both
additions in Eq.~\eqref{eq:boundary-correction} are finite power sums
with integer multiplicities. This together with Eq.~\eqref{eq:orig_sum} proves the first statement of the proposition asserted in Eq.~\eqref{eq:a-power-sum}.

It remains to prove Eq.~\eqref{eq:total-multiplicity}.
Combining Eq.~\eqref{eq:a-power-sum} with
Eq.~\eqref{eq:fifth-subsequence} gives for $\mathcal A_{5r}(a)$
\[
\mathcal A_{5r}(a)=\sum_{i=1}^km_i(\gamma_i^5)^r
=\bigl((q^2\alpha)^5\bigr)^r+\bigl((q^2\beta)^5\bigr)^r
\]
for all $r \geq 1$.
Group bases with the same fifth power by writing
\[
M_\delta=
\sum_{\substack{i \in \{1,\dots,k\},\\\gamma_i^5=\delta}}m_i,
\]
so the previous equation becomes
\[
\sum_{\delta \in \Delta} M_\delta\delta^r=\bigl((q^2\alpha)^5\bigr)^r+\bigl((q^2\beta)^5\bigr)^r,
\]
for all $r\geq 1$ 
where $\Delta\subseteq \CC^\times$ is a finite set. 

Recall that distinct nonzero power sequences are linearly independent which follows from the invertibility of the Vandermonde matrix. Comparing coefficients
therefore gives $M_{(q^2\alpha)^5}=M_{(q^2\beta)^5}=1$ and  $M_\delta=0$ for all other $\delta$ (in the special case $\alpha=\beta$ we have of course $M_{(q^2\alpha)^5}=2$ and all other $M_\delta=0$).
Summing over the finitely many bases on either side yields
\[
\sum_{i=1}^km_i
=\sum_{\delta\in \Delta} M_\delta=2.
\]
\end{proof}

The subsequence has therefore determined the total
multiplicity of the entire sequence. We now use the identity of the Walsh transform of the Dobbertin function we obtained in Lemma~\ref{lem:apn} to relate the values of $\mathcal A_r$ at indices $r$ and $2r$.

\begin{proposition}\label{prop:descent-constraints}
For every $r\geq1$,
\begin{equation}\label{eq:a-doubling}
\mathcal A_r(a)^2-\mathcal A_{2r}(a)=2q^{5r}.
\end{equation}
\end{proposition}

\begin{proof}
If $5\nmid r$, then $5\nmid2r$ as well. Hence
Eq.~\eqref{eq:corrected-sequence}, Proposition~\ref{prop:split-iden},
and Lemma~\ref{lem:apn} give in this case
\[
\mathcal A_r(a)^2-\mathcal A_{2r}(a)
=S_r(a)^2+S_{2r}(a)
=W_{q^r}(a)^2+W_{q^{2r}}(a)
=2q^{5r}.
\]
If $5\mid r$, then, by the definitions given in Eq.~\eqref{eq:corrected-sequence}, Eq.~\eqref{eq:curve_count} and Eq.~\eqref{eq:C-exceptional}, we have 
$\mathcal A_r(a)=-q^{2r}C_r(a)$ and
$\mathcal A_{2r}(a)=-q^{4r}C_{2r}(a)$. Theorem~\ref{thm:two-roots} therefore gives
\[
\mathcal A_r(a)^2-\mathcal A_{2r}(a)
=q^{4r}\bigl(C_r(a)^2+C_{2r}(a)\bigr)
=2q^{5r}.
\]
\end{proof}

We now know that $\mathcal A_r(a)$ is a finite power sum with total
multiplicity two that satisfies Eq.~\eqref{eq:a-doubling}.
 The following lemma shows that the
doubling identity forces a representation akin to the one given in Theorem~\ref{thm:two-roots} for the elliptic curve. 

\begin{lemma}\label{lem:rigidity}
Let $Q>0$ be real. Suppose a real sequence $(b_r)_{r\geq1}$
has a representation
\[
b_r=\sum_{i=1}^k m_i\gamma_i^r,
\]
where $m_i \in \ZZ$, $\gamma_i \in \CC^\times$, $\sum_{i=1}^k m_i=2$. 
If
\begin{equation}\label{eq:rigidity-identity}
b_r^2-b_{2r}=2Q^r
\end{equation}
for all $r\geq 1$
then there exist $u,v\in\CC^\times$ such that
$b_r=u^r+v^r$ and  $uv=Q$.
\end{lemma}

\begin{proof}
Without loss of generality assume that all $m_i \neq 0$ and all $\gamma_i$ are distinct. 

By linear independence of distinct nonzero power sequences,
we can compare coefficient-wise. 

Since the $b_r$ are real, conjugate $\gamma_i$ must have equal multiplicities $m_i$.
Expanding Eq.~\eqref{eq:rigidity-identity} yields
\begin{equation} \label{eq:expanded}
    \sum_{i,j \in \{1,\dots,k\}}m_im_j(\gamma_i\gamma_j)^r
-\sum_{i=1}^k m_i(\gamma_i^2)^r
-2Q^r=\sum_{\gamma'} n_{\gamma'}(\gamma')^r=0
\end{equation}
for all $r\geq 1$, where the $\gamma'$ collect all combinations of $\gamma_i^2,\gamma_i\gamma_j$ and $Q$. Since this equation holds for all $r\geq 1$, we have that all $n_{\gamma'}=0$.
Set \(R=\max_i|\gamma_i|\). We consider the coefficient $n_{\gamma'}$ in this equation for $\gamma'=R^2$. Write the nonreal
conjugate pairs of absolute value $R$ as
\[
\{z_1,\overline{z_1}\},\ldots,\{z_t,\overline{z_t}\},
\]
and recall that $m_{z_j}=m_{\overline{z_j}}$. The contributions to $n_{\gamma'}$ are then $m_R$, $m_{-R}$ (both can possibly be $0$ of course), $m_{z_j}$ for $1 \leq j \leq t$ and, if $R=\sqrt Q$, an additional term of $-2$. So
\begin{equation}\label{eq:outer-coefficient}
0=n_{R^2}=m_R^2+m_{-R}^2+2\sum_{j=1}^t m_{z_j}^2-m_R-m_{-R}=m_R(m_R-1)+m_{-R}(m_{-R}-1)+2\sum_{j=1}^t m_{z_j}^2.
\end{equation}
with an additional term of $-2$ if $R=\sqrt{Q}$.
Suppose first that \(R<\sqrt Q\). By definition of $R$ we have \(|\gamma_i|\le R\). Therefore
$|\gamma_i\gamma_j|\le R^2<Q$ for all choices of $i,j$. So $n_Q=-2$, since the contributions in Eq.~\eqref{eq:expanded} of the terms $-2Q^r$ cannot be canceled, contradicting that all $n_{\gamma'}=0$.

Suppose now that $R>\sqrt Q$.
Then the right hand side of Eq.~\eqref{eq:outer-coefficient} has to vanish, so $m_{z_j}=0$ for all $j$ and $m_R,m_{-R} \in \{0,1\}$. 
Consider Eq.~\eqref{eq:expanded} again but now look at the coefficient $n_{-R^2}$. We get
\[0=n_{-R^2}=2m_Rm_{-R}.\]
So $m_R=m_{-R}=1$ is impossible. Clearly, $m_R=m_{-R}=0$ is impossible since $R=\max_i|\gamma_i|$. So the only possibility left is $m_R=1$, $m_{-R}=0$ or  $m_R=0$, $m_{-R}=1$. Set $z$ as $R$ or $-R$ depending on which of $m_R$ and $m_{-R}$ are $1$.
Set
\[
V_r=b_r-z^r-(Q/z)^r=\sum_{i=1}^k m_i\gamma_i^r-z^r-(Q/z)^r.
\]
We want to show that $V_r=0$ for all $r\geq 1$.
All surviving bases $\gamma_i$ in the power representation of $V_r$ have absolute value less than $R$, since $z^r$ cancels out the one remaining term with absolute value $R$ in $b_r$ and $Q/R<R$ by our assumption on $R$.
We have
\begin{align*}
V_{2r}&=b_{2r}-z^{2r}-(Q/z)^{2r}=b_r^2-2Q^r-z^{2r}-(Q/z)^{2r}\\
&=(V_r+z^r+(Q/z)^r)^2-2Q^r-z^{2r}-(Q/z)^{2r} \\
&= V_r^2+2(z^r+(Q/z)^r)V_r.
\end{align*}
Assume that $V_r$ is not identically $0$. Let $D<R$ be the largest absolute value among all surviving bases in the power representation.
The terms in $2z^rV_r$ arising from bases of absolute value $D$ in $V_r$
have absolute value $RD$.
Bases occurring in $V_r^2$ have absolute value at most $D^2$ and in $(Q/z)^rV_r$ at most $DQ/R$. Similarly, all bases occurring in $V_{2r}$ have absolute value at most $D^2$. Both are less than $RD$, so the terms in $2z^rV_r$ cannot be canceled, giving a contradiction. So $V_r=0$ for all $r\geq 1$, and
\[
b_r=z^r+(Q/z)^r
\]
as claimed for this case.

Lastly, suppose $R=\sqrt Q$.

We extend the definition of the $b_r$ in the natural way to the negative numbers, and we show that Eq.~\eqref{eq:rigidity-identity}  holds for all integer indices $r$, not just positive ones. Indeed, consider
\begin{align*}
b_{r}^2-b_{2r}-2Q^r=
\sum_{i,j}m_im_j(\gamma_i\gamma_j)^{r}
-\sum_i m_i(\gamma_i^2)^r-2Q^r=\sum_{\gamma'} n_{\gamma'}(\gamma')^r,
\end{align*}
where the $\gamma'$ collect all combinations of $\gamma_i^2,\gamma_i\gamma_j$ and $Q$. By Eq.~\eqref{eq:rigidity-identity}, this sum of powers evaluates to $0$ for all $r\geq 1$, which then means that all $n_{\gamma'}=0$, so the sum also evaluates to zero for all negative $r$, so Eq.~\eqref{eq:rigidity-identity} holds for all integers.

Replacing every $\gamma_i$ by $Q/\gamma_i$ replaces $b_r$ by
\[
\widetilde b_r=Q^r b_{-r},
\]

Eq.~\eqref{eq:rigidity-identity} is then also preserved for the sequence $\widetilde b_r$ since

\begin{align*}
\widetilde b_r^2-\widetilde b_{2r}&=(Q^r b_{-r})^2-Q^{2r}b_{-2r}=Q^{2r}(b_{-r}^2-b_{-2r})\\
&=Q^{2r}\cdot2Q^{-r}=2Q^r.
\end{align*}

Assume first that there is $\gamma_i$ with absolute value less than $\sqrt{Q}$. Then the new sequence $\widetilde b_r$ contains a base with absolute value $Q/|\gamma_i|>\sqrt{Q}$. Now our argument from the case $R>\sqrt{Q}$ can be applied to the sequence $\widetilde b_r$, showing that there exist $u,v \in \CC^\times$ such that 
$\widetilde b_r=u^r+v^r$ and  $uv=Q$. But then \[
b_r=Q^r\widetilde b_{-r}
=(Q/u)^r+(Q/v)^r,
\]
and $(Q/u)(Q/v)=Q$, finishing the proof for this subcase. 

It thus only remains to consider the case where all $\gamma_i$ in the power representation of $b_r$ have absolute value exactly $\sqrt{Q}$. In this case, Eq.~\eqref{eq:outer-coefficient} becomes (recalling the additional term of $-2$) 
\begin{equation} \label{eq:sqrtq}
    m_R(m_R-1)+m_{-R}(m_{-R}-1)
+2\sum_{j=1}^t m_{z_j}^2=2.
\end{equation}
Assume first that there is one nonreal, conjugate pair $z_j,\overline{z_j}$ in this equation with $m_{z_j} \neq 0$. Observe that $m_R(m_R-1)$ cannot be negative, the same holds for the corresponding term containing $m_{-R}$. So the only way to satisfy the equation is to have a unique $m_{z_j}\in \{-1,1\}$ and $m_R,m_{-R} \in \{0,1\}$. Since all $\gamma_i$ have absolute value $R$, this means $b_r=m_R R^r+m_{-R}(-R)^r+m_{z_j}(z_j^r+\overline{z_j}^r)$. 
Since $\sum_{i=1}^k m_i=2$ this means $m_R+m_{-R}+2m_{z_j}=2$ yielding the only possibility $m_{z_j}=1$ and $m_R=m_{-R}=0$, so $b_r=z_j^r+\overline{z_j}^r$ with $z_j\overline{z_j}=R^2=Q$ as claimed.

The last case that remains is that there is no nonreal, conjugate pair in Eq.~\eqref{eq:sqrtq}, i.e., $\sum_{j=1}^t m_{z_j}^2=0$. Then $b_r=m_R R^r+m_{-R}(-R)^{r}$, and $\sum_{i=1}^k m_i=2$ implies immediately $m_R+m_{-R}=2$. Eq.~\eqref{eq:sqrtq} then simplifies to
\begin{align*}
2
&=m_R(m_R-1)+m_{-R}(m_{-R}-1)\\
&=(m_R+m_{-R})^2-(m_R+m_{-R})-2m_Rm_{-R}\\
&=2-2m_Rm_{-R}.
\end{align*}
So $m_Rm_{-R}=0$ and the only solutions are $m_R=2$, $m_{-R}=0$ and $m_R=0$, $m_{-R}=2$, so $b_r=2\gamma^r=\gamma^r+\gamma^r$ with $\gamma \in \{-R,R\}$, so $\gamma^2=R^2=Q$, finishing the proof.
\end{proof}

Applying Lemma~\ref{lem:rigidity} to
Propositions~\ref{prop:corrected-power-sum}
and~\ref{prop:descent-constraints}:
\begin{corollary} \label{cor:ar}
    We have for all $a\in L^\times$,  that there exist $u,v \in \CC^\times$  with $uv=q^5$ such that
    \begin{equation*}
\mathcal A_r(a)=u^r+v^r
\end{equation*}
for all $r\geq 1$.
\end{corollary}
In particular, Corollary~\ref{cor:ar} gives $\mathcal A_1(a)=-S_1(a)=-W_q(a)=u+v$, so the Walsh transform of the Dobbertin APN function is uniquely determined (for a fixed $a$) by $u$ and $v$.
Together with Eq.~\eqref{eq:fifth-subsequence}, we have
\[
\mathcal A_{5r}(a)=u^{5r}+v^{5r}=(q^2\alpha)^{5r}+(q^2\beta)^{5r}.
\]
Here $\alpha,\beta \in \CC^\times$ are the coefficients originating from the elliptic curve in Theorem~\ref{thm:two-roots}.
Note that 
$(q^2\alpha)(q^2\beta)=q^5$.
We will now show that indeed $\{u,v\}=\{q^2\alpha,q^2\beta\}$. We have to start by connecting $\mathcal A_{5r}$ to $\mathcal A_r$. We will do this using a Dickson polynomial.

More precisely, we will use the fifth Dickson polynomial, which can be directly defined as
\[
D_5(X,A)=X^5-5AX^3+5A^2X.
\]
It follows by direct expansion or using the well known properties of the Dickson polynomials that
\begin{equation*}
D_5(u+v,uv)=u^5+v^5.
\end{equation*}
This means that (using Corollary~\ref{cor:ar}) 
\begin{equation} \label{eq:dickson}
    D_5(\mathcal A_1(a),q^5)=\mathcal A_5(a)=(q^2\alpha)^{5}+(q^2\beta)^{5}.
\end{equation}

The following injectivity property shows that we can determine $\mathcal A_1(a)$ uniquely from that.

\begin{lemma}\label{lem:D5-injective}
For every $A\in\QQ$, the map $x\mapsto D_5(x,A)$ is injective
on $\QQ$.
\end{lemma}

\begin{proof}
Suppose $x\ne y$ are rational and
$D_5(x,A)=D_5(y,A)$. Put
\[
U=x^2+y^2-3A,\qquad V=xy+A.
\]
Expansion gives
\[
0=\frac{D_5(x,A)-D_5(y,A)}{x-y}
=U^2+UV-V^2.
\]
Consequently,
\[
(2U+V)^2=5V^2.
\]
Since $\sqrt5$ is irrational, $U=V=0$. But then
\[
(x-y)^2-5(x+y)^2=-4(U+3V)=0.
\]
The same irrationality forces $x-y=x+y=0$, contrary to $x\ne y$.
\end{proof}
This now allows us finally to relate $S_1(a)=W_q(a)$ to the elliptic curve, and thus to prove Theorem~\ref{thm:main}
\begin{theorem}[Copy of Theorem~\ref{thm:main}]
For $a\neq 1$, recall the elliptic curve given in Eq.~\eqref{eq:curve} \begin{equation*}
\mathcal E_a:\qquad Y^2+(\Tr_{L/K}(a)X+N_{L/K}(a))Y=X^3+e_2(a)X^2+\omega(a).
\end{equation*}
Then for $a\in L^\times$, we have
\[
W_q(a)= \begin{cases}
    q^2(\#\mathcal E_a(\F_q)-q-1), & a\neq 1, \\
    (-1)^{m+1}q^2(q+1), & a=1.
\end{cases}
\] 
\end{theorem}

\begin{proof}
By Eq.~\eqref{eq:dickson}, Corollary~\ref{cor:ar}, and
Eq.~\eqref{eq:fifth-subsequence},
\[
D_5(\mathcal A_1(a),q^5)
=(q^2\alpha)^5+(q^2\beta)^5
=D_5\bigl(q^2(\alpha+\beta),q^5\bigr)
=D_5(-q^2C_1(a),q^5).
\]
Both $\mathcal A_1(a)$ and $-q^2C_1(a)$ are integers.
The injectivity of the Dickson polynomial just proven in Lemma~\ref{lem:D5-injective} therefore gives
\begin{equation*}
\mathcal A_1(a)=-q^2C_1(a).
\end{equation*}
Since $\mathcal A_1(a)=-S_1(a)=-W_q(a)$ by the definition in Eq.~\eqref{eq:corrected-sequence}, we conclude
that $S_1(a)=q^2C_1(a)$. For $a=1$ this immediately gives $W_q(1)=(-1)^{m+1}q^2(q+1)$ using the special value from Eq.~\eqref{eq:C-exceptional}. In all other cases, we get the connection to the elliptic curve via Eq.~\eqref{eq:curve_count}.
\end{proof}

\section{Counting the points in the elliptic curve}
With Theorem~\ref{thm:main} proven, it remains to determine $\#\mathcal E_a(\F_q)$ for the elliptic curve defined by Eq.~\eqref{eq:curve}. This is quite standard.

For $a\notin\{0,1\}$, we now evaluate
$\#\mathcal E_a(K)$.

We first determine when the curve is supersingular and when it is ordinary. 
\begin{lemma}\label{lem:ordinary-supersingular}
Let $a\in L\setminus\{0,1\}$. The elliptic curve $\mathcal E_a$
is ordinary if and only if $\Tr_{L/K}(a)\ne0$, and supersingular if and only
if $\Tr_{L/K}(a)=0$.
\end{lemma}
\begin{proof}
The homogeneous equation of $\mathcal E_a$ is
\[
Y^2Z+\Tr_{L/K}(a)XYZ+N_{L/K}(a)YZ^2=X^3+e_2(a)X^2Z+\omega(a) Z^3.
\]
In characteristic two, a nonsingular plane cubic is supersingular
if and only if its $XYZ$ coefficient is zero
\cite[p.~332]{Hartshorn}.
Here that coefficient is $\Tr_{L/K}(a)$, proving the assertion.
\end{proof}

For $k\in K$, define the Kloosterman sum
\[
\Kl_K(k)=\sum_{u\in K^\times}\psi_K(u+k/u).
\]
We use the following result for the point counting of elliptic curves
\cite[Proposition~1 and Section~2, p.~881]{LachaudWolfmann1987}.
\begin{theorem}\label{thm:ordinary-point-count}
Let $K$ be a finite field of characteristic $2$. Every ordinary elliptic curve $E$ over $K$ is $K$-isomorphic to
a curve
\[
E_{d,k}:\qquad y^2+xy=x^3+dx^2+k,
\qquad d\in K,\quad k\in K^\times.
\]
In this case, we have
\begin{equation*}
\#E(K)=q+1+\psi_K(d)\Kl_K(k).
\end{equation*}
\end{theorem}
With this, we can prove Theorem~\ref{thm:main_explicit}.
\begin{proof}[Proof of Theorem~\ref{thm:main_explicit}]
For $a=0$, 
we have $W_q(0)=0$
For $a=1$, Theorem~\ref{thm:main} gives
\[
W_q(1)=(-1)^{m+1}q^2(q+1).
\]
So assume $a\in L\setminus\{0,1\}$.
If $\Tr_{L/K}(a)=0$, then $\omega(a)=N_{L/K}(a)e_3(a)$ and $N_{L/K}(a)\ne0$.
Substitute $X=u^2+e_2(a)$, and
with $z=Y/N_{L/K}(a)$ the curve equation becomes
\[
z^2+z
=\left(\frac{u^3+e_2(a)u}{N_{L/K}(a)}\right)^2+\frac{e_3(a)}{N_{L/K}(a)}.
\]
Since $z^2+z=t$ has $1+\psi_K(t)$ solutions and
$\psi_K(t^2)=\psi_K(t)$, each $u\in K$ contributes
\[
1+\psi_K\!\left(\frac{u^3+e_2(a)u+e_3(a)}{N_{L/K}(a)}\right)
\]
affine points. Summing over $u$ and adding $O$ proves
the stated formula for this case.

Now suppose $\Tr_{L/K}(a)\ne0$. Lemma~\ref{lem:auxiliary} gives
$\Tr_{L/K}(a)^2+N_{L/K}(a)/a\ne0$, so $k(a)\ne0$.
Expanding the norm gives
\[
N_{L/K}(\Tr_{L/K}(a)^2+N_{L/K}(a)/a)
=N_{L/K}(a)^4+N_{L/K}(a)^3\Tr_{L/K}(a)^3+e_2(a)N_{L/K}(a)^2\Tr_{L/K}(a)^4+\omega(a) \Tr_{L/K}(a)^6.
\]
The change of coordinates
\[
X=\Tr_{L/K}(a)^2x+\frac{N_{L/K}(a)}{\Tr_{L/K}(a)},
\qquad
Y=\Tr_{L/K}(a)^3y+\frac{N_{L/K}(a)^2}{\Tr_{L/K}(a)^3}
\]
therefore transforms the equation of $\mathcal E_a$ into
\[
y^2+xy=x^3+d(a)x^2+k(a).
\]
Theorem~\ref{thm:ordinary-point-count} then concludes the proof of this case.

Finally, if $a\in K\setminus\{0,1\}$, then $e_2(a)=0$, $\Tr_{L/K}(a)=a$,
and $N_{L/K}(a)=a^5$. Hence
$d(a)=a^2$ and 
$k(a)=\left(\frac{(a+1)^5}{a}\right)^2$, yielding the formula for this case.
\end{proof}

We now give an explicit evaluation of the cubic sums using
Carlitz's formulas
\cite[Theorems~1 and~2]{carlitz1979explicit}.
\begin{proposition}\label{prop:cubic-evaluation}
Let $a\in L\setminus\{0,1\}$ with $\Tr_{L/K}(a)=0$, and put
\[
A(a)=\sum_{u\in K}
\psi_K\!\left(\frac{u^3+e_2(a)u+e_3(a)}{N_{L/K}(a)}\right).
\]

\textup{(i)} Suppose $m$ is odd. Let $\rho\in K^\times$ be the
unique element with $\rho^3=N_{L/K}(a)$, and put $B=e_2(a)/\rho^2$.
If $\Tr_{K/\FF_2}(B)=0$, then $A(a)=0$.
If $\Tr_{K/\FF_2}(B)=1$, the equation
\[
z^4+z=B^2+1
\]
has a solution in $K$, and for any such solution,
\[
A(a)=(-1)^{(m^2-1)/8}\,2^{(m+1)/2}
\psi_K\!\left(\frac{e_3(a)}{N_{L/K}(a)}+z^3+z\right).
\]

\textup{(ii)} Suppose $m$ is even. If
\[
z^4+N_{L/K}(a)z=e_2(a)^2
\]
has no solution in $K$, then $A(a)=0$.
Otherwise, for any solution $z\in K$,
\[
A(a)=\psi_K\!\left(\frac{z^3+e_3(a)}{N_{L/K}(a)}\right)
\begin{cases}
(-1)^{m/2+1}2^{m/2+1},
  & N_{L/K}(a)\in(K^\times)^3,\\
(-1)^{m/2}2^{m/2},
  & N_{L/K}(a)\notin(K^\times)^3.
\end{cases}
\]
In particular,
\[
\#\mathcal E_a(K)=q+1+A(a).
\]
\end{proposition}

\begin{proof}
Extracting the constant term gives
\[
A(a)=\psi_K(e_3(a)/N_{L/K}(a))
\sum_{u\in K}\psi_K\!\left(N_{L/K}(a)^{-1}u^3+(e_2(a)/N_{L/K}(a))u\right).
\]
Apply \cite[Theorems~6.2 and~6.3, pp.~864 and~866]{Hou2007}
with $\alpha=1$, cubic coefficient $\kappa=N_{L/K}(a)^{-1}$, linear
coefficient $e_2(a)/N_{L/K}(a)$, and extension degree $m$.
For odd $m$, the substitution $u=\rho x$ gives $x^3+Bx$.
Square invariance yields
\[
\sum_{x\in K}\psi_K(x^3+Bx)
=
\sum_{x\in K}\psi_K(x^3+B^2x).
\]
For $B\ne0$, apply
\cite[Theorem~2, p.~15]{carlitz1979explicit}, keeping in mind the
identity $\left(\frac{2}{m}\right)=(-1)^{(m^2-1)/8}$ for Jacobi symbols.
This gives \textup{(i)}.
For $B=0$, the sum vanishes because cubing permutes $K$.

For even $m$, apply
\cite[Theorem~1, pp.~10--11]{carlitz1979explicit}
giving \textup{(ii)}.
\end{proof}

Using the well known characterization of the spectrum of Kloosterman sums, we can determine the spectrum of the Dobbertin APN function (Theorem~\ref{thm:walsh-inclusion}).
We start with a lemma
\begin{lemma}\label{lem:quadratic-surjectivity}
Put $H=\{x \in L \colon \Tr_{L/K}(x)=0\}$. Define
\[
\Phi:H\longrightarrow K^2,\qquad
\Phi(\theta)=
\bigl(e_2(\theta),\Tr_{L/K}(\theta^5)\bigr).
\]
Then $\Phi$ is surjective.
\end{lemma}

\begin{proof}
Recalling $e_1(\theta)=\Tr_{L/K}(\theta)$ and $e_5(\theta)=N_{L/K}(\theta)$, Newton's identities for elementary symmetric polynomials (or a direct elementary calculation) yield for $\theta \in H$
\begin{equation} \label{eq:newton}
    N_{L/K}(\theta)
=\Tr_{L/K}(\theta^5)
+e_2(\theta)\Tr_{L/K}(\theta^3).
\end{equation}

Consequently, $\Phi^{-1}(0,0)=\{0\}$.

Viewed over $\FF_2$, $H$ has dimension $4m$ and $K^2$ has dimension $2m$. Note that $\Phi$ is quadratic, since $e_2$ is quadratic and $\theta^5=\theta\theta^4$
is also quadratic over $\FF_2$.

Write its $2m$ coordinate functions as $f_1,\ldots,f_{2m}$.
For any $y\in\FF_2^{2m}$, the preimage size satisfies
\[
\begin{aligned}
\#\Phi^{-1}(y)
&\equiv
\sum_{x\in\FF_2^{4m}}
\prod_{i=1}^{2m}(1+f_i(x)+y_i)\\
&\equiv
\sum_{x\in\FF_2^{4m}}
\prod_{i=1}^{2m}(1+f_i(x))
\equiv\#\Phi^{-1}(0)
\pmod2.
\end{aligned}
\]
The second congruence holds because every term involving
some $y_i$ has degree at most $4m-2$ since all $f_i$ are quadratic, and a polynomial of
degree less than $4m$ sums to zero over $\FF_2^{4m}$.
Since the preimage of $0$ has one element, every preimage has odd
cardinality, and is in particular not empty.
\end{proof}

\begin{proof}[Proof of Theorem~\ref{thm:walsh-inclusion}]
First, $a=1$ gives $|W_q(1)|=q^2(q+1)$.

Let us first consider the set $\{q^2j \colon j \in \mathbb{N} \text{ odd}, |j|\leq 2\sqrt{q}\}$. Fix $k\in K^\times$.
Lemma~\ref{lem:quadratic-surjectivity} supplies $\theta\in H=\{x \in L \colon \Tr_{L/K}(x)=0\}$
such that
$e_2(\theta)=0,\Tr_{L/K}(\theta^5)=k.$
By Eq.~\eqref{eq:newton}, we then have $N_{L/K}(\theta)=k$.
Put
\[
a=\lambda(\theta+1),\text{ where }
\lambda^2=N_{L/K}(\theta+1)^{-1}.
\]
Since $\theta\ne0$ and $K\cap H=\{0\}$, this is well defined
and $a\notin K$. Write $\theta_i=\theta^{q^i}$ for $0\leq i<5$.
We have 
\[
\Tr_{L/K}(a)
=\lambda\bigl(\Tr_{L/K}(\theta)+\Tr_{L/K}(1)\bigr)
=\lambda(0+1)=\lambda,
\]
\[
\begin{aligned}
e_2(a)
&=\lambda^2\sum_{i<j}(\theta_i+1)(\theta_j+1)\\
&=\lambda^2\left(
e_2(\theta)+4\sum_{i=0}^4\theta_i+\binom52
\right)
=0,
\end{aligned}
\]
\[
N_{L/K}(a)
=\lambda^5N_{L/K}(\theta+1)
=\lambda^5\lambda^{-2}
=\lambda^3.
\]
Substituting into the definitions in
Theorem~\ref{thm:main_explicit} gives $d(a)=1$ and $k(a)=k$. Hence
\[
|W_q(a)|=q^2|\Kl_K(k(a))|.
\]
The range of Kloosterman sums is for $m\geq 2$ exactly \[
\{\Kl_K(k):k\in K^\times\}
=
\{j\in\mathbb Z:j\equiv-1\pmod4,\ |j|\le2\sqrt q\},
\]
see~\cite[Theorem~3.4]{LachaudWolfmann1990}, giving the stated values (the $m=1$ case can be verified by hand).

Let us now focus on the values in $\mathcal{U}_m$, and fix $\eta\in K^\times$
and $B\in K$.
Lemma~\ref{lem:quadratic-surjectivity} supplies a nonzero
$\theta\in H$ with $e_2(\theta)=B$: Indeed, 
use the target $(B,0)$ if $B\ne0$, and $(0,1)$ if $B=0$.
Set
\[
a=\lambda\theta, \text{where} 
\lambda^2=\bigl(\eta N_{L/K}(\theta)\bigr)^{-1}.
\]
Then $\Tr_{L/K}(a)=0$.
The substitution $u=\lambda x$ in the cubic-sum formula gives
\[
W_q(a)
=q^2\psi_K\bigl(\eta e_3(\theta)\bigr)
\sum_{x\in K}\psi_K\bigl(\eta(x^3+Bx)\bigr).
\]
Consequently,
\[
|W_q(a)|
=q^2\left|
\sum_{x\in K}\psi_K\bigl(\eta(x^3+Bx)\bigr)
\right|.
\]
The result then follows from the Walsh spectra of the Gold APN function, see Table~\ref{tab:apn-power-spectra}.
\end{proof}

\section{Sidon sets and two-error-correcting codes}

A subset $S$ of $\FF_2^t$ is Sidon if no four distinct elements sum to zero.
A result due to Czerwinski and Pott~\cite{CzerwinskiPott2026} constructs large Sidon sets in $\F_{2}^{2m-1}$ as well as binary linear codes with good parameters from APN functions in dimension $m$ with high linearity. Theorem~\ref{thm:main} in particular determines the linearity of the Dobbertin APN function, leading to the following result.
\begin{theorem}\label{thm:sidon-codes}
For every $m\geq1$, put
\[
M_m=2^{5m-1}+2^{3m-1}+2^{2m-1}.
\]
There exists a Sidon set of size $M_m$ in $\FF_2^{10m-1}$,
and there exist binary linear codes with parameters
\[
[M_m-1,M_m-10m,5]_2.
\]
\end{theorem}
\begin{proof}
By Theorem~\ref{thm:main}, the largest absolute Walsh
coefficient of the Dobbertin APN function is $\ell=q^2(q+1)$.
By \cite[Corollary~3.6]{CzerwinskiPott2026}, a Sidon set
can be constructed from an APN function by intersecting its graph with an affine hyperplane, the size of this Sidon set is
\[
\frac{2^{5m}+\ell}{2}=M_m.
\]
The intersection affinely spans the hyperplane: otherwise,
an additional independent affine equation would imply
$M_m\leq\ell$, contradicting
$M_m=(2^{5m}+\ell)/2>\ell$.

Let $S\subseteq\FF_2^{10m-1}$ be the resulting Sidon set.
Since $\binom{M_m}{3}>2^{10m-1}$, two distinct three-element
subsets of $S$ have the same sum. These subsets are disjoint:
if they shared exactly one point, the remaining four distinct
points would sum to zero, contradicting the Sidon property;
if they shared two points, equality of their sums would force
the subsets to coincide.

Thus there are six distinct points $s_1,\ldots,s_6\in S$ with
\[
s_1+\cdots+s_6=0.
\]
Choose $s_6$ as the translation origin and put
\[
S'=S+s_6=\{s+s_6:s\in S\}.
\]
Then $S'$ is Sidon, contains zero, and linearly spans
$\FF_2^{10m-1}$ because $S$ affinely spans that space.
Moreover, the five distinct nonzero points
$s_1+s_6,\ldots,s_5+s_6$ satisfy
\[
\sum_{i=1}^{5}(s_i+s_6)=0.
\]

We now follow the construction of Czerwinski and Pott~\cite[Section~5, before Theorem~5.1]{CzerwinskiPott2026}. Form a matrix $H$ whose columns are the elements of
$S'\setminus\{0\}$
and define
\[
C=\{c\in\FF_2^{M_m-1}:Hc^{\mathsf T}=0\}.
\]
The matrix $H$ has rank $10m-1$, so $C$ has length $M_m-1$
and dimension
\[
(M_m-1)-(10m-1)=M_m-10m.
\]
Its columns are distinct and nonzero, excluding codewords
of weights one and two. A dependence among three columns
would, together with zero, contradict the Sidon property
of $S'$, and a dependence among four columns would
contradict that property directly. Hence $C$ has minimum
distance at least five.

Finally, the columns corresponding to
$s_1+s_6,\ldots,s_5+s_6$ sum to zero and therefore give
a codeword of weight five, so the code has minimum distance exactly $5$.
\end{proof}

In $\FF_2^{10m-1}$, the classical constructions give Sidon sets
of size $2^{5m-1}+2^{\lfloor5m/2\rfloor}$, see
\cite[Remark~1.3]{CzerwinskiPott2026}.
When $m$ is odd, the inverse function on $\FF_{2^{5m}}$ is APN,
and its construction gives Sidon sets of size
\[
M_m^{\mathrm{inv}}
=2^{5m-1}+2\left\lfloor\frac{2^{5m/2+1}+1}{4}\right\rfloor
=2^{5m-1}+2^{5m/2}+O(1)
\]
by \cite[Theorem~4.3]{CzerwinskiPott2026}.
The subsequent improvement of the inverse parameters by one point,
noted in \cite[Section~6]{CzerwinskiPott2026}, has the same
asymptotic secondary term. These constructions have size
$2^{5m-1}+\Theta(2^{5m/2})$ in the dimensions where they apply.
In contrast, Theorem~\ref{thm:sidon-codes} gives
\[
M_m=2^{5m-1}+2^{3m-1}+2^{2m-1}
=2^{5m-1}+2^{3m-1}+o(2^{3m}).
\]
To the best of our knowledge, this gives the largest known
infinite family of Sidon sets in these ambient dimensions
asymptotically, with the best known secondary term.

The corresponding distance-five codes have length $N_m=M_m-1$
and dimension
\[
k_m=M_m-10m
=2^{5m-1}+2^{3m-1}+2^{2m-1}-10m.
\]
For the same ambient dimension $10m-1$, a spanning Sidon set
of size $M$ gives a code of dimension $M-10m$
\cite[Section~5]{CzerwinskiPott2026}.
Thus the improvement in Sidon-set size gives the same improvement
in code dimension: the secondary term grows from order
$2^{5m/2}$ to $2^{3m-1}$.

To place this in context, let $k_{\max}(N,5)$ be the largest
dimension of a binary linear code of length $N$ and minimum
distance at least five. Shortened double-error-correcting BCH
codes and the Hamming bound give, respectively,
\[
N-2\left\lceil\log_2(N+1)\right\rceil
\leq k_{\max}(N,5)
\leq N-\log_2\!\left(1+N+\binom N2\right),
\]
see \cite[Ch.~7, Theorem~3, p.~193; Ch.~1, Theorem~6,
p.~19]{MacWilliamsSloane1977}. Consequently,
\[
k_{\max}(N,5)=N-2\log_2 N+O(1).
\]
In fact, we can prove that the codes constructed through this construction are at most one dimension away from an optimal code for given ambient dimension and minimum distance $5$, which covers exactly the important codes that are $2$-error correcting.
\begin{theorem}\label{thm:code-near-optimal}
Let $m\geq1$, and put
\[
M_m=2^{5m-1}+2^{3m-1}+2^{2m-1},\qquad
N_m=M_m-1,\qquad k_m=M_m-10m.
\]
The code constructed in Theorem~\ref{thm:sidon-codes} has
dimension $k_m$ in $\FF_2^{N_m}$.
Every binary linear code $C\subseteq\FF_2^{N_m}$ with minimum
distance at least five satisfies
\[
\dim C\leq k_m+1.
\]
Thus our code is within one dimension of optimal among codes
in the same ambient space with minimum distance at least five.
\end{theorem}

\begin{proof}
Since $N_m>2^{5m-1}$, we have
\[
1+N_m+\binom{N_m}{2}
=\frac{N_m^2+N_m+2}{2}
>2^{10m-3}.
\]
The Hamming bound for codes~\cite[Ch.~1, Theorem~6, p.~19]{MacWilliamsSloane1977}
therefore yields
\[
2^{\dim C}\left(1+N_m+\binom{N_m}{2}\right)
\leq 2^{N_m},
\]
and hence
\[
\dim C<N_m-10m+3=k_m+2.
\]
As $\dim C$ is an integer, $\dim C\leq k_m+1$.
\end{proof}

Determining larger Sidon sets or more/better codes with minimum distance $5$ remains an open problem, a possible avenue is the construction of APN functions with linearity larger than the one of the Dobbertin APN function.
\section*{Declaration of AI use}
I obtained the main results of this paper without AI use in the first half of 2026, using the same rough approach shown here but also employing some more complicated tools, using function fields and the theory of zeta functions of elliptic curves. The proof shown here was streamlined and simplified using ChatGPT 6, in particular the use of Dwork's theorem in its current form was a suggestion made by the AI.

The same AI was also used for editorial help when writing the paper, and to provide literature references.
\section*{Acknowledgments}
I would like to thank Faruk G\"olo\u{g}lu, who suggested a connection to Kloosterman sums in the very early stages of this project, leading me to the idea to connect the sum to elliptic curves. I also thank Anne Canteaut, L\'{e}o Perrin, and Jules Baudrin for some discussions about this project during my visit to INRIA in Paris in the summer of 2026. These discussions in particular led to Lemma~\ref{lem:apn}, which streamlined the proof considerably. I thank INRIA and my hosts for their hospitality during my stay.
\bibliographystyle{amsplain}
\bibliography{bibliography} 
\end{document}